\documentstyle{amsppt}

\magnification=1200

\hsize = 6 true in
\vsize = 8.7 true in

\define\RR {\Bbb R}
\define\pa {\partial}
\define\ZZ {\Bbb Z}
\define\CC {\Bbb C}
\define\NN {\Bbb N}
\define\sgn {\operatorname{sgn}}

\define\bL {\bold{L}}
\define\F {\Cal{F}}
\define\Vol {\operatorname{Vol}}
\define\p {\partial}
\define\ep {\varepsilon}
\define\bW {\partial\Omega}
\define\Ombar {\overline{\Omega}}
\define\Ubar {\overline{U}}
\define\Ric {\operatorname{Ric}}
\define\bmo {\operatorname{bmo}}
\define\CO {\Cal{O}}
\define\OCO {\overline{\Cal{O}}}
\define\Lip {\operatorname{Lip}}
\define\End {\operatorname{End}}
\define\Mbar {\overline{M}}
\define\Tr {\operatorname{Tr}}
\define\CV {\Cal{V}}
\define\Exp {\operatorname{Exp}}
\define\wti {\widetilde{i}}
\define\dist {\operatorname{dist}}

\NoBlackBoxes

\centerline{\bf Boundary Regularity for the Ricci Equation,}
\vskip 7pt
\centerline{\bf Geometric Convergence, and Gel'fand's Inverse Boundary 
Problem}

$$\text{}$$

\centerline{\smc Michael Anderson,\quad Atsushi Katsuda,\quad
Yaroslav Kurylev,} 
\vskip 7pt
\centerline{\smc Matti Lassas,\quad and Michael Taylor}

$$\text{}$$

\centerline{\bf Abstract}
$\text{}$ \newline
This paper explores and ties together three themes.  The first is to 
establish regularity of a metric tensor, on a manifold with boundary,
on which there are given Ricci curvature bounds, on the manifold and 
its boundary, and a Lipschitz bound on the mean curvature of the boundary.
The second is to establish geometric convergence of a (sub)sequence of
manifolds with boundary with such geometrical bounds and also an upper
bound on the diameter and a lower bound on injectivity and boundary 
injectivity radius, making use of the first part.
The third theme involves the uniqueness and conditional stability
of an inverse problem proposed by Gel'fand, making essential
use of the results of the first two parts.

$$\text{}$$
{\bf 1. Introduction}
\newline {}\newline

The goals of this paper are to establish regularity, up to the boundary, 
of the metric tensor of a Riemannian manifold with boundary, under Ricci 
curvature bounds and control of the boundary's mean curvature; to apply this 
to results on Gromov compactness and geometric convergence in the category
of manifolds with boundary; and then to apply these results to the study of
an inverse boundary spectral problem introduced by I.~Gel'fand.

Regularity of the metric tensor away from the boundary has been studied and 
used in a number of papers, starting with [DTK].  The tack has been to 
construct local harmonic coordinates and use the fact that, in such harmonic
coordinates, the Ricci tensor has the form
$$
\Delta g_{\ell m}-B_{\ell m}(g,\nabla g)=-2 \Ric_{\ell m}.
\tag{1.0.1}
$$
Here $\Delta$ is the Laplace-Beltrami operator, applied {\it componentwise}
to the components of the metric tensor, and $B_{\ell m}$ is a quadratic 
form in $\nabla g_{ij}$, with coefficients that are smooth functions of
$g_{ij}$ as long as the metric tensor satisfies a bound $C_1|\eta|^2\le
g_{jk}(x)\eta^j\eta^k\le C_2|\eta|^2$, with $0<C_1\le C_2<\infty$.
If one is given information on the Ricci tensor, 
one can regard (1.0.1) as an elliptic PDE for the metric tensor,
and obtain information on its components, in harmonic coordinates.

The notion of compactness of a family of Riemannian manifolds and of geometric
convergence issues from work of J.~Cheeger [Ch] and M.~Gromov (cf.~[Gr],
the revised and translated version of his 1981 work).  The role of harmonic
coordinates in the study of such geometric convergence 
has been exploited in a number of papers.  It was used in [P] and in [GW] 
to obtain a compactness result, assuming a bound on the 
Riemann tensor, and some other geometric quantities.  In [An1] there was a 
successful treatment of compactness given a sup norm bound on the Ricci 
tensor, an upper bound on the diameter, and a lower bound on the injectivity 
radius, for a family of compact Riemannian manifolds of a fixed dimension.
Convergence was shown to hold, for a subsequence, in the $C^r$-topology,
for any $r<2$.  (A definition of geometric convergence is recalled in \S{3}.)

One of our motivations to extend the scope of these results to the category 
of manifolds with boundary arises naturally in the study of 
a class of inverse problems.  In these problems, 
one wants to determine the coefficients of some partial differential equation
in a bounded region via measurements of solutions to the PDE at the boundary.
Such problems arise in various areas, including geophysics, medical imaging,
and nondestructive testing.  One problem, formulated by I.~Gel'fand [Ge], 
consists of finding the shape of a compact manifold $\Mbar$ with boundary
$\pa M$ and the metric tensor on it from the spectral data on $\pa M$.
Namely, if $R_\lambda$ is the resolvent of the Neumann Laplacian $\Delta^N$ 
on $M$, the Gel'fand data consists of the restriction of the 
integral kernel $R_\lambda(x,y)$ of the resolvent to $x,y\in \pa M$,
as $\lambda$ varies over the resolvent set of $\Delta^N$.
Another formulation of Gel'fand's inverse problem will be given in \S{4}.

For such an inverse problem, the first issue to investigate is uniqueness.
In the context of $C^\infty$ metric tensors, this was established for the
Gel'fand problem in [BK1]
taking into account the unique continuation in [Ta].
See also [Bz], [NSU], [Be1] and [Nv] for the isotropic inverse problems.
As we will explain below it is important
to obtain uniqueness with much less regular coefficients. 

Once uniqueness results have been obtained, one has to face up to the
issue of {\it ill posedness} of the inverse problem.  That is, one can make
large changes in $\Mbar$ that have only small effects on boundary data 
obtained from examining the boundary behavior of the resolvent kernel 
mentioned above.  For example, given $(\Mbar,g)$, 
one could take an auxiliary manifold $X$, without boundary,
of the same dimension as $\Mbar$, 
remove a small ball from $X$ and from the interior of 
$\Mbar$, and connect these manifolds by a thin tube.  One is faced with the
task of {\it stabilizing} this ill posed inverse problem.  One ingredient 
in this process involves having some {\it a priori} knowledge of the 
quantities one is trying to determine, typically expressed in terms of
a priori bounds on these quantities in certain norms.  

An early result in this direction for the Gel'fand problem was given in [Al]
by G.~Alessandrini, who obtained {\it conditional stability} for the
operator $\text{div}\, \ep\, \text{grad}$ in a bounded domain in Euclidean 
space, where $\ep$ is a positive function (scalar conductivity), assumed to 
be bounded in some Sobolev space $H^s(M)$, with $s>0$.  See also [StU]
for a related result for an anisotropic metric tensor close to Euclidean.
Despite these successes, there is a clear need for coordinate-invariant 
constraints.

In the case of trying to determine an unknown 
Riemannian manifold with boundary $\Mbar$, from boundary spectral data, 
it is natural to make a priori hypotheses on geometrical properties of
$\Mbar$.  Furthermore, if one must make such a priori hypotheses, it 
is desirable to get by with as weak a set of hypotheses as possible.
There is then a tension between the desire to make weak a priori hypotheses 
and the need to establish uniqueness results.
(For preliminary results in this direction see [K2L], [Ka].)

Here we impose a priori sup norm bounds on the Ricci tensor of $\Mbar$,
and of $\pa M$.  This, together with a Lipschitz norm bound on the mean 
curvature of $\pa M\hookrightarrow \Mbar$, is shown in \S{2} to imply certain
regularity, up to the boundary, of the metric tensor of $\Mbar$, when one 
is in ``boundary harmonic coordinates'' (defined in \S{2}).  To be precise,
we obtain regularity in the Zygmund space $C^2_*(\Mbar)$, a degree of 
regularity better than $C^r$ for any $r<2$ and just slightly worse than 
$C^2$.  This result has the following important advantage over a $C^{2-\ep}$
estimate.  The Hamiltonian vector field associated with the metric tensor
has components with a log-Lipschitz modulus of continuity.  Hence, by
Osgood's theorem, it generates a uniquely defined geodesic flow, on the
interior of $\Mbar$, and also for geodesics issuing transversally from  
$\pa M$.  This property will be very important in \S{4}.
(We note that in the context of differential geometry Zygmund-type 
spaces go back to the habilitation thesis of B.~Riemann.)

In \S{3} we obtain a compactness result for families of
compact Riemannian manifolds, of dimension $n$, 
with boundary, for which there are fixed bounds on 
the sup norms of $\Ric_M$ and $\Ric_{\pa M}$, 
on the Lipschitz norm of the mean curvature of $\pa M$, and on the diameter, 
and fixed lower bounds on the injectivity and boundary 
injectivity radius.  We show that a sequence of such Riemannian manifolds
has a subsequence, converging in the $C^r$-topology, for all $r<2$,
whose limit $(\Mbar,g)$ has metric tensor in $C^2_*(\Mbar)$.

In \S{4} we study Gel'fand's inverse boundary problem, recast in the form of 
an inverse boundary spectral problem.  We show that, having boundary spectral
data, we can recognize whether a given function $h \in C(\pa M)$ has the 
form $h(z)=r_x(z)=\text{dist}(x,z)$, for some $x\in\Mbar$, all 
$z\in\pa M$, thus recovering the image in $C(\Mbar)$ of $\Mbar$
under the boundary distance representation.  
Such a representation, whose use was initiated in [Ku] and [KuL],
plays an important role in the uniqueness proof, but for it to work
we need to know that geodesics from points in $\pa M$, pointing normal to
the boundary, are uniquely defined.  As noted above, this holds when the 
metric tensor in in $C^2_*(\Mbar)$, and we obtain a uniqueness result in this 
category.  This fits in perfectly with the compactness result of \S{3}, to 
yield a result on stabilization of this inverse problem.

Section 5 is devoted to the proof of several elliptic regularity results, of 
an apparently non-standard nature, needed for some of the finer results of
\S{2}.

$\text{}$ \newline
{\smc Remark}.  A number of classes of function spaces arise naturally in
our analysis.  These include spaces $C^r(\Mbar)$, mentioned above.  Here,
if $r=k+\sigma,\ k\in\ZZ^+,\ \sigma\in (0,1)$, $C^r(\Mbar)$ consists of 
functions whose derivatives of order $k$ satisfy a H{\"o}lder condition,
with exponent $\sigma$.  The Zygmund spaces $C^r_*(\Mbar)$ coincide with 
$C^r(\Mbar)$ for $r\in (0,\infty)\setminus\ZZ^+$, and form a complex 
interpolation scale.  We also encounter $L^p$-Sobolev spaces, $H^{s,p}(M)$
and Besov spaces $B^s_{p,p}(\pa M)$, and $\text{bmo}(M)$, the localized
space of functions of bounded mean oscillation.
Basic material on these spaces can be 
found in [Tr1] and Chapters 2, 3 of
[Tr2], in Chapter 13 of [T1], and in Chapter 1 of [T2].

$\text{}$ \newline
{\smc Acknowledgment}. Thanks to R.~Mazzeo for useful discussions, 
particularly regarding material in \S{2.1}.  Thanks also to T.~Sakai,
T.~Sunada, and Y.D.~Burago for discussions on the role of geometric 
convergence, and to M.~Gromov, G.~Uhlmann, and E.~Somersalo for interest
and support.  M.~Anderson's work was partially supported by NSF grant
DMS-0072591.  A.~Katsuda's work was partially supported by a Grant in Aid for
Scientific Research (C)(2) No.~14540081 of JSPS.  M.~Lassas's work was
partially funded by the Academy of Finland.
M.~Taylor's work was partially supported by NSF grant
DMS-0139726.  Part of this work was done at MSRI and at the Oberwolfach
RiP Program, and their support is gratefully acknowledged.

$$\text{}$$
{\bf 2. Boundary Regularity for the Ricci Equation}
\newline {}\newline

In this section we establish the key results on local regularity at
the boundary of a metric tensor on which there are Ricci curvature
bounds and a Lipschitz bound on the mean curvature.  Our set-up
is the following.

Let $\Cal{B}$ be a ball about $0\in\RR^n,\ \Omega=\Cal{B}\cap\{x:x^n>0\}.$
Let $\Sigma=\Cal{B}\cap \{x:x^n=0\}$ and set $\Ombar=\Omega\cup\Sigma$.
Let $g$ be a metric tensor on $\Ombar$, and denote by $h$ its restriction to
$\Sigma$.  We make the following hypotheses:
$$
\align
g_{jk}&\in H^{1,p}(\Omega),\ \text{ for some }\ p>n,
\tag{2.0.1} \\
h_{jk}&\in H^{1,2}(\Sigma),\quad 1\le j,k\le n-1,
\tag{2.0.2} \\
\Ric^{\Omega}&\in L^\infty(\Omega), \tag{2.0.3} \\
\Ric^{\Sigma}&\in L^\infty(\Sigma), \tag{2.0.4} \\
H&\in \Lip(\Sigma), \tag{2.0.5}
\endalign
$$

Here $H$ denotes the mean curvature of $\Sigma\subset \Ombar$, i.e.,
$H=\text{Tr}\ A/(n-1)$, where $A$ is the Weingarten map, a section of
$\End(T\Sigma)$.  Our goal is to establish the following result.

\proclaim{Theorem 2.1} Under the hypotheses (2.0.1)--(2.0.5), 
given $z\in\Sigma$, 
there exist local harmonic coordinates 
on a neighborhood $\Ubar$ of $z$ in $\Ombar$ with respect to which 
$$
g_{jk}\in C^2_*(\Ubar).
\tag{2.0.6}
$$
\endproclaim

Here $C^2_*(\Ubar)$ is a Zygmund space, as mentioned in \S{1}.
The harmonic coordinates for which (2.0.6) holds are arbitrary 
coordinates $(u^1,\dots,u^n)$ satisfying $\Delta u^j=0$ on a chart
not intersecting $\Sigma$.  On a neighborhood of a point in $\Sigma$,
these coordinates are ``boundary harmonic coordinates,'' which are 
defined as follows.  We require $(u^1,\dots,u^n)$ to be defined and 
regular of class at least $C^1$ on a neighborhood of $z$ in $\Ombar$,
and $\Delta u^j=0$.  We require that $v^j=u^j|_{\Sigma}$ be harmonic 
on $\Sigma$, i.e., annihilated by the Laplace-Beltrami operator of 
$\Sigma$ with its induced metric tensor.  We require $u^n$ to vanish
on $\Sigma$, and we require $(u^1,\dots,u^n)$ to map a neighborhood
of $z$ in $\Ombar$ diffeomorphically onto $\Ombar$.

Let us note that the hypotheses (2.0.1)--(2.0.2) imply that various 
curvature tensors are well defined.  If $(g_{jk})$ is the $n\times n$
matrix representation of a metric tensor in a coordinate system, $(g^{jk})$
its matrix inverse, the connection 1-form $\Gamma$ is given by
$$
\Gamma^a{}_{bj}=\frac{1}{2} g^{am}(\pa_j g_{bm}+\pa_b g_{jm}-\pa_m g_{bj}).
\tag{2.0.7}
$$
The Riemann tensor is then given by
$$
\Cal{R}=d\Gamma+\Gamma\wedge\Gamma.
\tag{2.0.8}
$$
It is a matrix valued 2-form with components $R^a{}_{bjk}$.  We see that
$$
\aligned
g_{jk}\in C(\Ombar)\cap H^{1,2}(\Omega)&\Longrightarrow 
\Gamma\in L^2(\Omega),\ R^a{}_{bjk}\in H^{-1,2}(\Omega)+L^1(\Omega) \\
&\Longrightarrow \Ric_{bk}\in H^{-1,2}(\Omega)+L^1(\Omega).
\endaligned
\tag{2.0.9}
$$
The hypothesis (2.0.1) is stronger than the hypothesis in (2.0.9).  
It implies $g_{jk}\in C^r(\Ombar)$ for some $r>0$, so (2.0.9) is applicable
both to $g_{jk}$ on $\Omega$ and,
in view of (2.0.2), to $h_{jk}$ on $\Sigma$.  Furthermore,
$$
g_{jk}\in H^{1,p}(\Omega),\ p>n\Longrightarrow 
R^a{}_{bjk},\ \Ric_{bk},\ \Ric^j{}_k,\ S\in H^{-1,p}(\Omega),
\tag{2.0.10}
$$
where $\Ric^j{}_k=g^{jb}\Ric_{bk}$ and $S=\Ric^j{}_j$ is the scalar curvature
of $\Omega$.  We mention parenthetically that one can use the fact that 
pointwise multiplication gives a map
$$
H^{1,2}\times H^{-1,2}\longrightarrow H^{-1,p'},\quad 
\forall\ p'<\frac{n}{n-1},
\tag{2.0.11}
$$
to obtain
$$
g_{jk}\in C(\Ombar)\cap H^{1,2}(\Omega)\Longrightarrow
\Ric^j{}_k,\ S\in H^{-1,p'}(\Omega).
\tag{2.0.12}
$$
However, we will not make use of (2.0.12) here.

We next consider the implication of (2.0.1) for the Weingarten map 
associated to $\Sigma\hookrightarrow \Ombar$.  The unit normal $N$ to
$\Sigma$ is a vector field with coefficients
$$
N^j=\frac{1}{\sqrt{g^{nn}}}\, g^{jn}\bigr|_{\Sigma},
\tag{2.0.13}
$$
which by the trace theorem belongs to the Besov space $B^{1-1/p}_{p,p}
(\Sigma)$.  It follows that the Weingarten map has the property
$$
A\in B^{-1/p}_{p,p}(\Sigma),
\tag{2.0.14}
$$
as a consequence of (2.0.1).  Thus we have a priori that $H\in B^{-1/p}_{p,p}
(\Sigma)$, and the hypothesis (2.0.5) strengthens this condition on $H$,
in a fashion that is natural for the desired conclusion of Theorem 2.1.

Our approach to the proof of Theorem 2.1 is to obtain the result as a
regularity result for an elliptic boundary problem.  We use the PDE
(1.0.1) (the ``Ricci equation'') for the components of the metric tensor,
in boundary harmonic coordinates, and use Dirichlet boundary conditions on 
some components of $g_{jk}$ and Neumann boundary conditions on complementary
components; see (2.1.8) and (2.1.16)--(2.1.17) for a more precise
description.

We will approach the proof of Theorem 2.1 in stages.  In \S{2.1} we prove
that the conclusion (2.0.6) holds when the hypothesis (2.0.1) is
strengthened to $g_{jk}\in C^{1+s}(\Ombar)$, for some $s>0$.  In \S{2.2}
we replace (2.0.1) by the hypothesis that $g_{jk}\in H^{1,p}(\Omega)$ for 
some $p>2n$.  In \S{2.3} we prove the full strength version of Theorem 2.1.
These stages serve to isolate three rather different types
of arguments, each of which is needed to prove 
Theorem 2.1, but which are perhaps more digestible when presented 
separately.  Section 2.4 has some complementary results on the degree
of regularity of the harmonic coordinates mentioned in Theorem 2.1.
In section 2.5 we demonstrate the non-branching of geodesics for 
metric tensors satisfying (2.0.6), including geodesics starting at a 
boundary point, in a direction transversal to the boundary.  We also
discuss examples of branching geodesics, for metric tensors only mildly
less regular than those of (2.0.6), extending some examples of [Ha].

The version of Theorem 2.1 established in \S{2.1} is already useful for
the results of \S\S{3--4}, and the reader particularly interested in 
\S\S{3--4} could skip \S\S{2.2--2.3}, on first reading.  However, the 
hypothesis (2.0.1) has a ``natural'' quality that we believe makes the 
additional effort required to work with it worthwhile.  The arguments in 
\S\S{2.2--2.3} require several elliptic regularity results that do not
seem to be standard in the literature that we know, and their proofs are
collected later, in \S{5}.

$\text{}$

\heading
\S{2.1}: First regularity result
\endheading

$\text{}$

Here we prove that the conclusion of Theorem 2.1 holds when the hypotheses 
(2.0.1)--(2.0.2) are strengthened a bit.

\proclaim{Proposition 2.1.1} In the setting of Theorem 2.1, replace
hypotheses (2.0.1)--(2.0.2) by 
$$
g_{jk}\in C^{1+s}(\Ombar),\ \text{ for some }\ s\in (0,1),
\tag{2.1.1} 
$$
and retain hypotheses (2.0.3)--(2.0.5).  Then the conclusion (2.0.6) holds.
\endproclaim

To begin our demonstration,
let $h$ denote the metric tensor induced on $\Sigma$, with respect to which 
(2.0.4) holds.  By (2.1.1), $h_{jk}\in C^{1+s}(\Sigma)$, so there exist local 
harmonic coordinates $v^1,\dots,v^{n-1}$ on a neighborhood $\Cal{O}$ of $z$ 
in $\Sigma$.  Now we can find harmonic $u^1,\dots,u^{n-1}$ on a neighborhood 
of $z$ in $\Ombar$ such that $u^j=v^j$ on $\Cal{O}$.  Also we can find 
$u^n$, harmonic in $\Ombar$, with $u^n|_{\Sigma}=0$ and 
arrange that $\pa_{x^n}u^n(z)\neq 0$.  We will have
$$
u^1,\dots,u^n\in C^{2+s}(\Ubar).
\tag{2.1.2}
$$
We have $du^1(z),\dots,du^n(z)$ linearly independent, so, after perhaps 
further shrinking $\Ubar$ we have a harmonic coordinate chart on $\Ubar$, 
a set we relabel as $\Ombar$.  As mentioned below the statement of 
Theorem 2.1, this is what we call a set 
of {\it boundary harmonic coordinates}.  In these new coordinates, (2.1.1) 
and (2.0.3)--(2.0.5) are preserved.

Now in harmonic coordinates the metric tensor satisfies the elliptic PDE
$$
\Delta g_{\ell m}=F_{\ell m},
\tag{2.1.3}
$$
where $\Delta$ acts componentwise on $g_{\ell m}$, as
$$
\Delta u=g^{-1/2}\pa_j(g^{1/2} g^{jk}\pa_ku),\quad g=\det (g_{jk}),
\tag{2.1.4}
$$
and
$$
F_{\ell m}=B_{\ell m}(g,\nabla g)-2\Ric^\Omega_{\ell m}.
\tag{2.1.5}
$$
Here $B_{\ell m}$ is a quadratic form in $\nabla g$ with coefficients that
are rational functions of $g_{jk}$.  Thus, from (2.1.1) and (2.0.3) we have
$$
F_{\ell m}\in L^\infty(\Omega),
\tag{2.1.6}
$$
and the coefficients of $\Delta$ have the same degree of regularity as 
$g_{jk}$ in (2.1.1).

Now, if $j,k\le n-1$, then well known local regularity results on $\Sigma$
following from (2.0.4) give
$$
g_{jk}\bigr|_{\Sigma}=h_{jk}\in H^{2,p}(\Sigma),\quad \forall\ p<\infty,
\tag{2.1.7}
$$
but in fact there is the following refinement of (2.1.7), established in
Proposition III.10.2 of [T2]:
$$
g_{jk}\bigr|_{\Sigma}=h_{jk}\in \frak{h}^{2,\infty},\quad 1\le j,k\le n-1.
\tag{2.1.8}
$$
Here $\frak{h}^{2,\infty}$ denotes the bmo-Sobolev space of functions whose
derivatives of order $\le 2$ belong to $\bmo$, the localized space of 
functions of bounded mean oscillation.
We establish the following (after perhaps shrinking $\Ombar$ to a smaller
neighborhood of $z$).

\proclaim{Lemma 2.1.2} Under our working hypotheses we have, in the harmonic 
coordinate system $(u^1,\dots,u^n)$,
$$
g_{jk}\in C^2_*(\Ombar),\quad 1\le j,k\le n-1.
\tag{2.1.9}
$$
\endproclaim
\demo{Proof} First, extend $F_{\ell m}$ by $0$ on $\Cal{B}\setminus\Omega$
and solve $\widetilde{\Delta}w_{\ell m}=F_{\ell m}$ on a neighborhood of $0$
in $\Cal{B}$, where we obtain $\widetilde{\Delta}$ in the form (2.1.4) with 
$g^{jk}$ extended across $\Sigma$ to $\tilde{g}^{jk}\in C^{1+s}(\Cal{B}).$
Local elliptic regularity results imply 
$$
w_{\ell m}\in\frak{h}^{2,\infty}(\Cal{B})\subset C^2_*(\Cal{B}).
\tag{2.1.10}
$$
It follows that $w_{\ell m}|_{\Sigma}\in C^2_*(\Sigma)$ and, via (2.1.8),
$$
g_{jk}-w_{jk}\bigr|_{\Sigma}=b_{jk}\in C^2_*(\Sigma),\quad
j,k\le n-1,
\tag{2.1.11}
$$
while
$$
\Delta(g_{jk}-w_{jk})=0\ \text{ on }\ \Omega.
\tag{2.1.12}
$$
Given our assumed regularity of the coefficients of $\Delta$, standard
Schauder results give
$$
b_{jk}\in C^r(\Sigma)\Longrightarrow g_{jk}-w_{jk}\in C^r(\Ombar),\quad
1<r<2,\ 2<r<2+s.
\tag{2.1.13}
$$
Actually, the case $1<r<2$ is perhaps not so classical, but see [Mo1], 
Theorem 7.3 or [GT], Corollaries 8.35--8.36.
From here, an interpolation argument gives
$$
b_{jk}\in C^2_*(\Sigma)\Longrightarrow g_{jk}-w_{jk}\in C^2_*(\Ombar).
\tag{2.1.14}
$$
See [T1], Chapter 13, \S{8}, particularly (8.37), for interpolation in this 
context.  This establishes (2.1.9).
\enddemo

To continue, following [An2], we switch over to PDE for $g^{\ell m}$.
Parallel to (2.1.3), we have
$$
\Delta g^{\ell m}=B^{\ell m}(g,\nabla g)+2(\Ric^\Omega)^{\ell m}=F^{\ell m},
\tag{2.1.15}
$$
and (2.1.1) and (2.0.3) give $F^{\ell m}\in L^\infty(\Omega)$.  
We take $m=n$ and proceed to derive Neumann-type boundary
conditions for the components $g^{\ell n},\ 1\le\ell\le n$.  
In fact, as we will show,
$$
N g^{nn}=-2(n-1)Hg^{nn},\quad \text{on }\ \Sigma,
\tag{2.1.16}
$$
and, for $1\le\ell\le n-1$,
$$
N g^{\ell n}=-(n-1)Hg^{\ell n}+\frac{1}{2}
\frac{1}{\sqrt{g^{nn}}}\, g^{\ell k}\pa_k g^{nn},\quad \text{on }\ \Sigma.
\tag{2.1.17}
$$
Here $H$ is the mean curvature of $\Sigma$, which we assume satisfies (2.0.5),
and $N$ is the unit normal field to $\Sigma$, pointing inside $\Omega$.

To compute (2.1.16)--(2.1.17), we use
$$
g^{\ell m}=\langle \nabla u^\ell,\nabla u^m\rangle,\quad
N=\frac{\nabla u^n}{|\nabla u^n|}=\frac{1}{\sqrt{g^{nn}}}\nabla u^n,
\tag{2.1.18}
$$
and
$$
Ng^{\ell n}=\langle \nabla_N\nabla u^\ell,\nabla u^n\rangle
+\langle \nabla u^\ell,\nabla_N \nabla u^n\rangle.
\tag{2.1.19}
$$
We also use the fact that $u^\ell$ is harmonic on $\Omega$ and $u^\ell
\bigr|_{\Sigma}=v^\ell$ is harmonic on $\Sigma$ ($0$ if $\ell=n$).

Note that if $\{e_j:1\le j\le n-1\}$ is an orthonormal frame on $\Sigma$
and $X$ a vector field on $\Ombar$ (say both having coefficients in 
$C^1(\Ombar)$) then
$$
\text{div}\ X\bigr|_{\Sigma}=\sum\limits_{j=1}^{n-1} \langle\nabla_{e_j}
X,e_j\rangle +\langle \nabla_NX,N\rangle.
\tag{2.1.20}
$$
In particular, for $X_\ell=\nabla u^\ell$, we have $\text{div}\, X_\ell
=\Delta u^\ell=0$, so the first term on the right side of (2.1.19) is equal 
to $-\sqrt{g^{nn}}$ times
$$
\sum\limits_{j=1}^{n-1} \langle \nabla_{e_j}X_\ell,e_j\rangle.  
\tag{2.1.21}
$$
Let us set
$$
X_\ell=X_\ell^N+X_\ell^T,\quad X_\ell^N=\langle X_\ell,N\rangle N=
\varphi N,\quad X^T_\ell=\nabla v^\ell,
\tag{2.1.22}
$$
with $X^T_\ell$ tangent to $\bW$ and $\varphi = g^{ln}/ \sqrt{g^{nn}}$.
Since $\sum_j \langle \nabla_{e_j}\nabla v^\ell,e_j\rangle=\text{div}\, 
\nabla v^\ell=\Delta v^\ell=0$, we have (2.1.21) equal to
$$
\aligned
\sum\limits_{j=1}^{n-1} \langle\nabla_{e_j}(\varphi N),e_j\rangle
&=\varphi\sum\limits_j \langle\nabla_{e_j}N,e_j\rangle \\
&=\varphi \sum \langle Ae_j,e_j\rangle \\ 
&=(n-1)H\, \frac{g^{\ell n}}{\sqrt{g^{nn}}},
\endaligned
\tag{2.1.23}
$$
so the first term on the right side of (2.1.19) is equal to 
$-(n-1)Hg^{\ell n}$.  The case $\ell=n$ gives (2.1.16),
since the two summands in (2.1.19) are then the same.

To continue when $\ell\neq n$, we note that
$$
\langle \nabla u^\ell,\nabla_N \nabla u^n\rangle
=\langle N,\nabla_{X_\ell}\, \nabla u^n\rangle,
\tag{2.1.24}
$$
with $X_\ell=\nabla u^\ell$.  In fact, generally a 1-form $\eta$ satisfies
$$
d\eta(X,Y)=\langle \nabla_X\eta,Y\rangle-\langle\nabla_Y\eta,X\rangle,
$$
and applying this to $\eta=du^n$
and $X=X_\ell, \, Y=X_n$ gives (2.1.24).  Now
$$
X_\ell g^{nn}=2\langle \nabla_{X_\ell}\nabla u^n,\nabla u^n\rangle,
$$
so (2.1.24) is equal to 
$$
\frac{1}{2\sqrt{g^{nn}}}\, \langle X_\ell,\nabla g^{nn}\rangle
=\frac{1}{2\sqrt{g^{nn}}}\, g^{\ell k}\pa_k g^{nn},
\tag{2.1.25}
$$
which gives (2.1.17).

Having (2.1.15)--(2.1.17), we can establish further regularity of the 
functions $g^{\ell n}$.  

\proclaim{Lemma 2.1.3} In the harmonic coordinate system $(u^1,\dots,u^n)$,
we have
$$
g^{\ell n}\in C^2_*(\Ombar),\quad 1\le\ell\le n.
\tag{2.1.26}
$$
\endproclaim
\demo{Proof} As in Lemma 2.1.2, we extend $g^{jk}$ to 
$\tilde{g}^{jk}\in C^{1+s}(\Cal{B})$ and extend the right side of (2.1.15) by
$0$ and produce a solution 
$$
w^{\ell n}\in \frak{h}^{2,\infty}(\Cal{B})\subset C^2_*(\Cal{B})
\tag{2.1.27}
$$ 
to $\Delta w^{\ell m}=F^{\ell m}$.
Then $g^{\ell n}-w^{\ell n}$ satisfies
$$
\Delta (g^{\ell n}-w^{\ell n})=0, 
\quad \text{on }\ \Omega,
\tag{2.1.28}
$$
and we have
$$
N(g^{nn}-w^{nn})\bigr|_{\Sigma}\in C^1_*(\Sigma),
\tag{2.1.29}
$$
by (2.1.16) and (2.1.27), plus the regularity $N\in C^{1+s}(\Sigma)$.
As in Lemma 2.1.2, we can apply an interpolation 
argument to Schauder-type estimates (see \S{5.3}) and get
$$
g^{nn}-w^{nn}\in C^2_*(\Ombar),
\tag{2.1.30}
$$
and hence (2.1.26) holds for $\ell=n$.  Having this, we get from (2.1.17) 
that
$$
N(g^{\ell n}-w^{\ell n})\bigr|_{\Sigma}\in
C^1_*(\Sigma),
\tag{2.1.31}
$$
which gives $g^{\ell n}-w^{\ell n}\in C^2_*(\Ombar)$, so (2.1.26) holds for
all $\ell$.
\enddemo

The final step is to verify that Lemmas 2.1.2 and 2.1.3 yield regularity of 
$g_{\ell n}$.  

\proclaim{Lemma 2.1.4} In the setting of Lemmas 2.1.2--2.1.3,
$$
g_{n\ell}=g_{\ell n}\in C^2_*(\Ombar).
\tag{2.1.32}
$$
\endproclaim
\demo{Proof}  Let $g=\det(g_{jk})$
and let $A_{\ell m}$ be the determinant of the $(n-1)\times (n-1)$ matrix
formed by omitting column $\ell$ and row $m$ from the matrix $(g_{jk})$.
Then
$$
g^{jk}=\frac{(-1)^{j+k}}{g} A_{jk}.
\tag{2.1.33}
$$
By Lemma 2.1.2, $A_{nn}\in C^2_*(\Ombar)$.  Applying Lemma 2.1.3
to $g^{nn}$ (which is $>0$) we have
$$
g=A_{nn}/g^{nn}\in C^2_*(\Ombar).
\tag{2.1.34}
$$
Then it follows that
$$
A_{\ell n}=A_{n\ell}= (-1)^{n+\ell} g\, g^{n\ell}\in C^2_*(\Ombar),\quad 1\le \ell\le n.
\tag{2.1.35}
$$
Another way of putting this is the following.  Let 
$$
h_{jk}= g_{jk}, \quad 1 \leq j,k \leq n-1; \quad h = \det(h_{jk})
\tag{2.1.36}
$$
and $(h^{jk})$ be the matrix inverse to $(h_{jk})$. 
Then $A_{\ell n}, \ell \leq n-1,$ can be written in the form 
$$
A_{\ell n} = (-1)^{n-1+\ell} g_{j n} \, h h^{j \ell}.
\tag{2.1.37}
$$
Now the regularity and positive-definiteness of $(h_{jk})_{1 \leq j,k 
\leq n-1}$
applied to (2.1.37) yield
$$
(g_{n1},\dots,g_{n\,n-1})\in C^2_*(\Ombar).
\tag{2.1.38}
$$
Finally, the identity
$$
g_{jn} g^{jn}=1,
\tag{2.1.39}
$$
the regularity of $g^{jn}$ in (2.1.26) and of $g_{jn}$ for $j\le n-1$ in
(2.1.38), plus the fact that $g^{nn}>0$, yield
$$
g_{nn}\in C^2_*(\Ombar),
\tag{2.1.40}
$$
proving the lemma, and completing the proof of Proposition 2.1.1.
\enddemo

$\text{}$

\heading
\S{2.2}: First improvement
\endheading

$\text{}$

In this section and the next we obtain regularity with a weaker a priori
hypothesis than (2.1.1).  As we noted above, the results of \S{2.1}
suffice for the applications in \S{3}, but these improvements are quite
natural (if not trivial to implement)
and surely have the potential for applications elsewhere.

Here we do strengthen the hypothesis (2.0.1) to some degree.  
Namely we assume:
$$
g_{jk}\in H^{1,p}(\Omega),\quad p>2n,
\tag{2.2.1}
$$
We retain hypothesis (2.0.2), i.e.,
$$
h_{jk}\in H^{1,2}(\Sigma),\quad 1\le j,k\le n-1.
\tag{2.2.2}
$$
Here $\Omega$ and $\Sigma$ are as in \S{1} and $h_{jk}=g_{jk}|_{\Sigma}$,
for $1\le j,k\le n-1$.  Note that (2.2.1) implies $g_{jk}\in C^r(\Ombar)$
with $r=1-n/p>0$, and hence $h_{jk}\in C^r(\Sigma)$.  

Our strategy is to show that the hypotheses (2.2.1)--(2.2.2) together with 
(2.0.3)--(2.0.5) imply $g_{jk}\in C^{1+s}(\Ombar)$ for some $s>0$,
so Proposition 2.1.1 applies.  In fact, we will expand the scope of the 
investigation here, and establish this conclusion under the following
hypotheses, which are weaker than (2.0.3)--(2.0.5):
$$
\align
\Ric^\Omega&\in L^{p_1}(\Omega),\quad p_1>n, \tag{2.2.3} \\
\Ric^\Sigma&\in L^{p_2}(\Sigma),\quad p_2>n-1, \tag{2.2.4} \\
H&\in C^s(\Sigma),\quad s>0. \tag{2.2.5}
\endalign
$$
Our goal in this section is to prove:

\proclaim{Proposition 2.2.1} Under the hypotheses (2.2.1)--(2.2.5), 
given $z\in\Sigma$,
there exist local harmonic coordinates on a neighborhood $\Ubar$ of $z$ 
in $\Ombar$ with respect to which
$$
g_{jk}\in C^{1+s}(\Ubar),
\tag{2.2.6}
$$
for some $s>0$.
\endproclaim

As before, we begin by constructing local harmonic coordinates $v^1,\dots,
v^{n-1}$ on a neighborhood $\Cal{O}$ of $z$ in $\Sigma$.  Knowing that
$h_{jk}\in C^r(\Sigma)$, we can do this, and making use also of hypothesis
(2.2.2) we have
$$
v^j\in C^{1+r}(\Cal{O})\cap H^{2,2}(\Cal{O}),
\tag{2.2.7}
$$
by Proposition 9.4 in Chapter III of [T2].  It follows that (2.2.2) persists
in this new coordinate system.  As a consequence of the fact that $h_{jk}\in
B^{1-1/p}_{p,p}(\Sigma)$, we also have
$$
v^j\in B^{2-1/p}_{p,p}(\Cal{O}).
\tag{2.2.8}
$$
This result is established in \S{5.1}.

Next we find harmonic functions $u^1,\dots,u^{n-1}$ on a neighborhood
$\Ubar$ of $z$ in $\Ombar$ such that $u^j=v^j$ on $\Cal{O}$ and we find $u^n$,
harmonic in $\Ombar$, with $u^n|_\Sigma=0$ and arrange that $\pa_{x^n}u^n(z)
\neq 0$.  Given (2.2.7)--(2.2.8) and the hypothesis (2.2.1), we claim that
$$
u^1,\dots,u^n\in C^{1+r}(\Ubar)\cap H^{2,p}(U).
\tag{2.2.9}
$$
The fact that $u^j\in C^{1+r}(\Ubar)$ follows from Corollaries 8.35--8.36
of [GT], or [Mo1], Theorem 7.3.
The fact that $u^j\in H^{2,p}(U)$ is established in \S{5.2}.
We have $du^1(z),\dots,$ $du^n(z)$ linearly independent, so 
after perhaps further
shrinking $\Ubar$ we have a harmonic coordinate chart on $\Ubar$, which we
relabel $\Ombar$.  In these new coordinates, (2.2.1)--(2.2.5) are preserved.

In fact, now that we have switched to harmonic coordinates, we can improve 
(2.2.2), making use of (2.2.4). 
It follows from Proposition 10.1 in Chapter III of [T2] that
$$
h_{jk}\in H^{2,p_2}(\Sigma),\quad 1\le j,k\le n-1.
\tag{2.2.10}
$$
In particular,
$$
h_{jk}\in C^{1+\sigma}(s),\ \text{ for some }\ s>0.
\tag{2.2.11}
$$
We may as well suppose $s \in (0,r)$.  Now we can prove:

\proclaim{Lemma 2.2.2} In the harmonic coordinate system $(u^1,\dots,u^n)$,
$$
g_{\ell m}\in C^{1+s}(\Ombar),\quad 1\le \ell,m\le n-1.
\tag{2.2.12}
$$
\endproclaim
\demo{Proof} We know $g_{\ell m}$ solves the Dirichlet problem
$$
\Delta g_{\ell m}=F_{\ell m},\quad g_{\ell m}\bigr|_\Sigma=h_{\ell m},
\quad 1\le \ell,m\le n-1,
\tag{2.2.13}
$$
where
$$
F_{\ell m}=B_{\ell m}(g,\nabla g)-2\Ric^\Omega_{\ell m}.
\tag{2.2.14}
$$
From (2.2.1) and (2.2.3) we have
$$
F_{\ell m}\in L^{q_1}(\Omega),\quad q_1=\min(p/2,p_1)>n,
\tag{2.2.15}
$$
and the coefficients $a^{jk}=g^{1/2}g^{jk}$ of $\Delta$ are known to belong 
to $H^{1,p}(\Omega)\subset C^r(\Ombar)$.

The next step in the proof is by now familiar.  Extend $F_{\ell m}$ by $0$
on $\Cal{B} \setminus \Omega$ and solve $\Delta v_{\ell m}=F_{\ell m}$ on a 
neighborhood $\Cal{V}$ of $z$ in $\Ombar$ 
with $v_{\ell m}\in C^{1+s}(\Cal{V}),\, s>0$.
Then $w_{\ell m}=g_{\ell m}-v_{\ell m}$ solves
$$
\Delta w_{\ell m}=0\ \text{ on }\ \Omega,\quad w_{\ell m}\bigr|_\Sigma=
h_{\ell m}-v_{\ell m}\bigr|_\Sigma\in C^{1+s}(\Sigma),
\tag{2.2.16}
$$
and the previously cited results of [Mo1] and [GT] yield $w_{\ell m}\in
C^{1+s}(\Ombar)$, hence (2.2.12).
\enddemo

It remains to show that
$$
g_{\ell n}=g_{n\ell}\in C^{1+s}(\Ombar),\quad 1\le \ell\le n,
\tag{2.2.17}
$$
with $s>0$.  In fact, if we show that
$$
g^{\ell n}=g^{n\ell}\in C^{1+s}(\Ombar),\quad 1\le \ell\le n,
\tag{2.2.18}
$$
then an argument parallel to the proof of Lemma 2.1.4 yields (2.2.17).

As before, we have
$$
\Delta g^{\ell n}=B^{\ell n}(g,\nabla g)+2(\Ric^\Omega)^{\ell n}=F^{\ell n}.
\tag{2.2.19}
$$
As in (2.2.15), we have
$$
F^{\ell n}\in L^{q_1}(\Omega),\quad q_1=\min(p/2,p_1)>n.
\tag{2.2.20}
$$
However, this time it is not so straightforward to produce the Neumann-type
boundary conditions (2.1.16)--(2.1.17).  

Consider (2.1.16).  The right side is well defined; 
we have $Hg^{nn}|_{\Sigma}\in C^s(\Sigma)$, for some $s>0$.  
As for $N$, the unit normal field to 
$\Sigma$ is also H{\"o}lder continuous of class $C^r$. 
But applying $N$ to $g^{nn}\in H^{1,p}(\Omega)$
does not yield an object that can be evaluated on $\Sigma$.  One has the 
same problem with the left side of (2.1.17), and the right side of (2.1.17)
is also problematic.

Our next goal is to show that a weak formulation of the Neumann boundary 
condition is applicable.  Generally, the weak formulation of 
$$
\Delta w=F,\quad Nw\bigr|_{\Sigma}=G
\tag{2.2.21}
$$
is that for all test functions $\psi$, i.e., all $\psi\in C^\infty(\Ombar)$
with compact support (intersecting $\Sigma$ but not the rest of $\pa\Omega$),
$$
\int\limits_{\Omega} \langle\nabla w,\nabla \psi\rangle\, dV=
-\int\limits_{\Omega} F\psi\, dV-\int\limits_{\Sigma} G\psi\, dS.
\tag{2.2.22}
$$
Here $dV$ is the volume element on $\Omega$ and $dS$ the area element on
$\Sigma$, both determined by the metric tensor in the usual fashion.  
Note that the left side of (2.2.22) is well defined for all test functions
$\psi$ whenever $\nabla w\in L^1(\Omega)$ and the right side of (2.2.22)
is well defined whenever $F\in L^1(\Omega)$ and $G\in L^1(\Sigma)$.

\proclaim{Lemma 2.2.3} The function $w=g^{nn}$ satisfies (2.2.22), with
$$
\aligned
F&=F^{nn}\in L^{q_1}(\Omega),\quad q_1>n, \\
G&=-2(n-1)Hg^{nn}\bigr|_{\Sigma}\in C^r(\Sigma),\quad r>0.
\endaligned
\tag{2.2.23}
$$
Hence the result that $g^{nn}\in H^{1,p}(\Omega)$ for some $p>2n$
is improved to
$$
g^{nn}\in C^{1+s}(\Ombar),\quad s>0.
\tag{2.2.24}
$$
\endproclaim

That (2.2.22) holds in this context seems a natural 
generalization of (2.1.16), but the proof, given below,
requires some work.  Given that (2.2.22) holds,
the regularity result (2.2.24) follows from the results \S{5.3}.

Once Lemma 2.2.3 is established, we see that the right side of (2.1.17) is 
well defined, and we can formulate:

\proclaim{Lemma 2.2.4} For $1\le\ell\le n-1$, the function $w=g^{\ell n}$
satisfies (2.2.22), with
$$
\aligned
F&=F^{\ell n}\in L^{q_1}(\Omega),\quad q_1>n, \\
G&=-(n-1)Hg^{\ell n}\bigr|_{\Sigma}+\frac{1}{2}\frac{1}{\sqrt{g^{nn}}}\,
g^{\ell k}\pa_k g^{nn}\bigr|_{\Sigma}\in C^r(\Sigma),\quad r>0.
\endaligned
\tag{2.2.25}
$$
Hence the result $g^{\ell n}\in H^{1,p}(\Omega)$, for some $p>2n$,
is improved to
$$
g^{\ell n}\in C^{1+s}(\Ombar),\quad s>0.
\tag{2.2.26}
$$
\endproclaim

To set up the proof of Lemmas 2.2.3--2.3.4, 
let $\Omega_c=\{x\in\Omega:u^n(x)>c\}\subset \Omega$, 
for small $c>0$, and let $\Sigma_c=\{x\in\Omega:u^n(x)=c\}$.  
Since $u^n\in C^{1+s}$ the surfaces $\Sigma_c$ are uniformly 
$C^{1+s}$-smooth.  Parallel to (2.2.10), we have from (2.2.3) that
$$
g_{jk}\in H^{2,p_1}_{\text{loc}}(\Omega),
\tag{2.2.27}
$$
so $Ng^{\ell n}|_{\Sigma_c}$ is well defined
for small $c>0$.  Calculations parallel to (2.1.18)--(2.1.23) give
$$
Ng^{nn}\bigr|_{\Sigma_c}=-2(n-1)H_c g^{nn},
\tag{2.2.28}
$$
where $H_c$ denotes the mean curvature of $\Sigma_c$, i.e., $(n-1)H_c=
\text{Tr}\, A_c$, where $A_c$ denotes the Weingarten map of $\Sigma_c$.  
Note that, for $X,Y$ tangent to $\Sigma_c$,
$$
\langle A_cX,Y\rangle=\langle \nabla_XN,Y\rangle,\quad 
N=\frac{\nabla u^n}{|\nabla u^n|},\quad (\nabla u^n)^j=g^{jn}.
\tag{2.2.29}
$$
The assumption (2.2.1) implies
$$
N\in H^{1,p}(\Omega),\quad N\bigr|_{\Sigma_c}\in B^{1-1/p}_{p,p}
(\Sigma_c),\quad A_c\in B^{-1/p}_{p,p}(\Sigma_c).
\tag{2.2.30}
$$
For a fixed $c>0$, from (2.2.27),
we have $N|_{\Sigma_c}\in B^{2-1/p_1}_{p_1,p_1}
(\Sigma_c)$, hence $A_c\in B^{1-1/p_1}_{p_1,p_1}(\Sigma_c)$, 
but (2.2.30) holds uniformly as $c\rightarrow 0$,
and if $\Sigma_c$ is identified with $\Sigma$
via the $(u^1,\cdots,u^{n-1})$-coordinates,
we have $A_c$ continuous in $c$ as $c\rightarrow 0$, in the space 
$B^{-1/p}_{p,p}(\Sigma)$.  

Now pick a test function $\psi$.  From (2.2.28) we have, for each (small)
$c>0$,
$$
\int\limits_{\Omega_c} \langle \nabla g^{nn},\nabla \psi\rangle \, dV
=-\int\limits_{\Omega_c} F^{nn}\psi\, dV
-\int\limits_{\Sigma_c} G_c \psi\, dS,
\tag{2.2.31}
$$
where  $F^{nn}, G_c$ are as in (2.2.23) with $\Sigma _c$ instead of $\Sigma$.
We let $c\rightarrow 0$.  Since we already
have $\nabla g^{nn}\in L^2(\Omega)$, the left side of (2.2.31) converges to
$$
\int\limits_{\Omega} \langle \nabla g^{nn},\nabla\psi\rangle\, dV,
\tag{2.2.32}
$$
and the first term on the right side of (2.2.31) converges to
$$
-\int\limits_{\Omega} F^{nn}\psi\, dV.
\tag{2.2.33}
$$
Finally, from (2.2.30) we have that $H_c\rightarrow H_0$ in $B^{-1/p}_{p,p}
(\Sigma)$, with $p>2n$, which, given our knowledge at this point that 
$g^{jk}\in H^{1,p}(\Omega)$, is more than enough to
imply that the last term in (2.2.31) converges to
$$
-\int\limits_{\Sigma} G\psi\, dS.
\tag{2.2.34}
$$
This proves Lemma 2.2.3.

To proceed, we have, for $1\le\ell\le n-1,\ c>0$, by a calculation parallel
to (2.1.18)--(2.1.25),
$$
\aligned
Ng^{\ell n}\bigr|_{\Sigma_c}=\ &-(n-1)H_cg^{\ell n}
-\sqrt{g^{nn}} \Delta_{\Sigma_c} u^\ell\bigr|_{\Sigma_c} \\
&+\frac{1}{2} \frac{1}{\sqrt{g^{nn}}}\, g^{\ell k} \pa_k g^{nn}
\bigr|_{\Sigma_c}.
\endaligned
\tag{2.2.35}
$$
Here $\Delta_{\Sigma_c}$ is the Laplace operator on the surface 
$\Sigma_c$, with its induced Riemannian metric tensor.
Hence, given a test function $\psi$, we have, for each $c>0$,
$$
\int\limits_{\Omega_c} \langle\nabla g^{\ell n},\nabla\psi\rangle\, dV
=-\int\limits_{\Omega_c} F^{\ell n}\psi\, dV - \int\limits_{\Sigma_c}
G\psi\, dS,
\tag{2.2.36}
$$
with $F^{\ell n}$ as in 
(2.2.25)
and $G_c$ given by the right side of (2.2.35).  Again we let 
$c\rightarrow 0$, and since we know $\nabla g^{\ell n}\in L^2(\Omega)$, the 
left side of (2.2.36) converges to
$$
\int\limits_{\Omega} \langle\nabla g^{\ell n},\nabla \psi\rangle\, dV.
\tag{2.2.37}
$$
Again the first term on the right side of (2.2.36) converges to
$$
-\int\limits_{\Omega} F^{\ell n} \psi\, dV.
\tag{2.2.38}
$$

The last term in (2.2.36) is equal to
$$
\aligned
&(n-1)\int\limits_{\Sigma_c} H_c g^{\ell n}\psi\, dS
+\int\limits_{\Sigma_c} \sqrt{g^{nn}} \bigl(\Delta_{\Sigma_c}u^\ell
\bigr|_{\Sigma_c}\bigr) \psi\, dS \\
&-\frac{1}{2} \int\limits_{\Sigma_c} \frac{1}{\sqrt{g^{nn}}}
g^{\ell k} \pa_k g^{nn} \psi\, dS.
\endaligned
\tag{2.2.39}
$$
As $c\rightarrow 0$, the first term in (2.2.39) converges to
$$
(n-1)\int\limits_{\Sigma} Hg^{\ell n}\psi\, dS,
\tag{2.2.40}
$$
by the same arguments as above (via (2.2.30)).  Next (making use of (2.2.1)) 
we have, for some $p>2n$,
$$
u^\ell\in H^{2,p}(\Omega),\quad u^\ell\bigr|_{\Sigma_c}
\in B^{2-1/p}_{p,p}(\Sigma_c),\quad \Delta_{\Sigma_c}
u^\ell\bigr|_{\Sigma_c}\in B^{-1/p}_{p,p}(\Sigma_c),
\tag{2.2.41}
$$
with uniform bounds and convergence as $c\rightarrow 0$, so the second term
in (2.2.39) converges to $0$ as $c\rightarrow 0$.

Finally, since we already have $g^{nn}\in C^{1+s}(\Ombar)$, plus $g^{\ell k}
\in H^{1,p}(\Omega)$, the convergence of the last term in (2.2.39) as 
$c\rightarrow 0$ follows, so Lemma 2.2.4 is proven.

$\text{}$

\heading
\S{2.3}: Proof of Theorem 2.1
\endheading

$\text{}$

In this section we finish the proof of Theorem 2.1
assuming that $n < p \leq 2n$. 
Going further, we extend Proposition 2.2.1 as follows.

\proclaim{Theorem 2.3.1} Replace hypothesis (2.2.1) by
$$
g_{jk}\in H^{1,p}(\Omega),\quad n<p\leq 2n,
\tag{2.3.1}
$$
and retain hypotheses (2.2.2)--(2.2.5).  Then, 
given $z\in\Sigma$, 
there exist
local harmonic coordinates on a neighborhood $\Ubar$ of $z$ in $\Ombar$
with respect to which
$$
g_{jk}\in C^{1+s}(\Ubar)
\tag{2.3.2}
$$
for some $s>0$.
\endproclaim

To begin, we have local harmonic coordinates $v^1,\dots,v^{n-1}$ on a 
neighborhood $\Cal{O}$ of $z$ in $\Sigma$ satisfying (2.2.7)--(2.2.8) and then 
local harmonic coordinates $u^1,\dots,u^n$ as in \S{2.2}, satisfying (2.2.9),
and in these new coordinates (2.3.1) and (2.2.2)--(2.2.5) are preserved.
We also continue to have (2.2.10)--(2.2.11).  We next establish the variant
of Lemma 2.2.2 that holds in this context.

\proclaim{Lemma 2.3.2} In the harmonic coordinate system $(u^1,\dots,u^n)$, 
we have
$$
g_{\ell m}\in H^{1,r_1}(\Omega),\quad 1\le \ell,m\le n-1,
\tag{2.3.3}
$$
with
$$
r_1=\frac{q_1}{1-q_1/n},\quad \text{for any}\, \,
 \frac{n}{2} <q_1 <\frac{p}{2}.
\tag{2.3.4}
$$
\endproclaim
\demo{Proof} We continue to have (2.2.13)--(2.2.15), except that in (2.2.15) 
no longer have  $q_1>n$ (rather $n/2<q_1<n$).  As for the 
regularity of $g_{\ell m}|_\Sigma$, we continue to have (2.2.11).  Extend 
$F_{\ell m}$ by $0$ on $\Cal{B}\setminus\Omega$ and solve $\Delta v_{\ell m}=
F_{\ell m}$ on a neighborhood $\Cal{V}$ of $z$ in $\Cal{B}$.  
This time we can 
say $v_{\ell m}\in H^{1,r_1}(\Cal{V})$, with $r_1$ as in (2.3.4).  Since
$L^{q_1}\subset H^{-1,r_1}$, this follows from Proposition 1.10 in [T2],
Chapter III.  Then $w_{\ell m}=g_{\ell m}-v_{\ell m}$ satisfies (2.2.16), so 
as before we have $w_{\ell m}\in H^{1,r_1}(\Ombar)$, and this gives (2.3.3).
Note that
$$
q_1>\frac{n}{2}\Longrightarrow r_1>2q_1.
\tag{2.3.5}
$$
\enddemo

To proceed, we note that the results (2.2.27)--(2.2.34) hold under our 
relaxed hypotheses on $p$, so (2.2.23) holds, except we have $q_1$ as in 
(2.3.4) instead of $q_1>n$.  With this, we can prove the following.

\proclaim{Lemma 2.3.3} For all $\ep>0$, we have
$$
g^{nn}\in H^{1,r_1-\ep}(\Omega).
\tag{2.3.6}
$$
\endproclaim
\demo{Proof} We first extend $F^{nn}$ by $0$ on
$\Cal{B}\setminus\Omega$ 
and solve
$$
\Delta g^{nn}_0=F^{nn},\quad g^{nn}_0\in H^{1,r_1}(\Cal{V}),
\tag{2.3.7}
$$
on a neighborhood $\Cal{V}$ of $z$ in $\Cal{B}$, as in 
the proof of Lemma 2.3.2.
In fact, by Proposition 5.2.2, we have
$$
g^{nn}_0\in H^{2,q_1}(\Cal{V}),
\tag{2.3.8}
$$
since $1/p+1/p\le 1/q_1$ in this case.  Of course (2.3.8) implies (2.3.7), 
but it also implies
$$
\nabla g^{nn}_0\bigr|_{\Sigma}\in B^{1-1/q_1}_{q_1,q_1}(\Sigma)\subset
L^{s_1-\ep}(\Sigma),\quad s_1=\frac{(n-1)q_1}{n-q_1},
\tag{2.3.9}
$$
the latter inclusion holding for all $\ep>0$.

Hence $g^{nn}=g^{nn}_0+g^{nn}_1$ where
$$
\Delta g^{nn}_1=0,\quad Ng^{nn}_1=G-Ng^{nn}_0.
\tag{2.3.10}
$$
Here $G=-2(n-1)Hg^{nn}|_{\Sigma}$, as in (2.2.23), so we know
$$
G\in C^r(\Sigma),\quad Ng^{nn}_0\in L^{s_1-\ep}(\Sigma),
\tag{2.3.11}
$$
for all $\ep>0$.  It follows from Proposition 5.5.2 that
$$
g^{nn}_1\in H^{1,r_1-\ep}(\Omega),\quad \forall\ \ep>0,
\tag{2.3.12}
$$
which gives (2.3.6).  For use below we also record the non-tangential maximal
function estimate
$$
(\nabla g^{nn}_1)^*\in L^{s_1-\ep}(\Sigma),
\tag{2.3.13}
$$
also established in \S{5.5}.  The meaning of the left side of (2.3.13) is
the following.  First, we have $\nabla g_1^{nn}$ continuous on the 
{\it interior} of $\Omega$.  Next, for $x\in\Sigma$,
$$
(\nabla g^{nn}_1)^*(x)=\sup\limits_{y\in\Gamma_x}\, |\nabla g^{nn}_1(y)|,
$$
where $\Gamma_x=\{y\in\Omega: d(y,x)\le 2\, d(y,\Sigma)\}$.
\enddemo

\proclaim{Lemma 2.3.4} For $1\le \ell\le n-1$, and for all $\ep>0$, we have
$$
g^{\ell n}\in H^{1,r_1-\ep}(\Omega).
\tag{2.3.14}
$$
\endproclaim
\demo{Proof} Using (2.3.8) and (2.3.13), we can extend the analysis in 
(2.2.35)--(2.2.41), to conclude that $g^{\ell n}$ is a weak solution to
$$
\Delta g^{\ell n}=F^{\ell n},\quad Ng^{\ell n}=G,
\tag{2.3.15}
$$
with $G$ as in (2.2.25), i.e.,
$$
G=-(n-1)Hg^{\ell n}\bigr|_{\Sigma}+\frac{1}{2} \frac{1}{\sqrt{g^{nn}}}
g^{\ell k}\pa_k g^{nn}\bigr|_{\Sigma}.
\tag{2.3.16}
$$
From what we know so far, we have
$$
G\in L^{s_1-\ep}(\Sigma),\quad \forall\ \ep>0.
\tag{2.3.17}
$$
As in the proof of Lemma 2.3.3, we can then write 
$$
g^{\ell n}=g^{\ell n}_0+g^{\ell n}_1,\quad g^{\ell n}_0\in H^{2,q_1}(\Omega),
\tag{2.3.18}
$$
with
$$
\Delta g^{\ell n}_1=0,\quad Ng^{\ell n}_1=G-Ng^{\ell n}_0\in L^{s_1-\ep}
(\Sigma),
\tag{2.3.19}
$$
which by an analysis parallel to that of (2.3.10)--(2.3.11) gives
$$
g^{\ell n}_1\in H^{1,r_1-\ep}(\Omega)
\tag{2.3.20}
$$
and proves (2.3.14).
\enddemo

Now an argument parallel to the proof of Lemma 2.1.4 gives
$$
g_{jk}\in H^{1,r_1-\ep}(\Omega),\quad \forall\ \ep>0,
\tag{2.3.21}
$$
for all $j,k\le n$.  Now, for $\ep$ small and $q_1$ close to $p/2$, 
$$
r_1-\ep>p,
\tag{2.3.22}
$$
an improvement over the hypothesis (2.3.1), as long as $p\le 2n$.  Thus 
replacing (2.3.1) by (2.3.22) and iterating this argument a finite number of
times, we establish that actually (2.2.1) holds.  

This proves Theorem 2.3.1.  It also reduces Theorem 2.1 to Proposition 2.1.1,
and hence proves Theorem 2.1.

$\text{}$ \newline
{\smc Remark}.  Theorem 2.1 remains valid if we change conditions
(2.0.3)--(2.0.5) into 
$$
\Ric ^{\Omega} \in \bmo (\Omega), \quad \Ric ^{\Sigma} \in  \bmo
(\Sigma),
\quad H \in C^1_*(\Sigma).
\tag{2.3.23}
$$

$\text{}$

\heading
\S{2.4}: Complements on coordinates
\endheading

$\text{}$

Under the hypotheses of Theorem 2.1, we know there are local boundary
harmonic coordinates with respect to which $g_{jk}\in C^2_*(\Ombar)$.
Here we show that, with respect to such coordinates, any {\it other}
boundary harmonic coordinates are smooth of class $C^3_*(\Ombar)$, so
$\Ombar$ has the structure of a $C^3_*$-manifold.  Since the interior
behavior is simpler to establish than the behavior at the boundary, we
confine our analysis to the following result.

\proclaim{Proposition 2.4.1} Assume the metric tensor on $\Ombar$ is of
class $C^2_*(\Ombar)$.  Let $u\in H^{1,2}(\Omega)$ solve
$$
\Delta u=f\in C^1_*(\Ombar),\quad u\bigr|_{\bW}=h\in C^3_*(\bW).
\tag{2.4.1}
$$
Then $u\in C^3_*(\Ombar)$.
\endproclaim
\demo{Proof} Classical Schauder estimates readily give $u\in 
C^{2+s}_*(\Ombar)$ for all $s<1$, as in (2.1.2); we just need to go a little 
further.  The key to success is to abandon our former practice (which worked
so well) of writing the Laplace-Beltrami operator $\Delta$ in divergence form,
and instead write it in non-divergence form:
$$
g^{jk}\pa_j\pa_k u=f_1,\quad u\bigr|_{\bW}=h,
\tag{2.4.2}
$$
where
$$
f_1=f-g^{-1/2}\pa_j(g^{1/2}g^{jk})(\pa_k u).
$$
The hypothesis on $g_{jk}$, plus the
current handle we have on $u$, gives $g^{-1/2}\pa_j(g^{1/2}g^{jk})$ $(\pa_ku)
\in C^{1}_*(\Ombar)$, hence $f_1\in C^1_*(\Ombar)$, under the hypothesis
on $f$ in (2.4.1).  Now the hypothesis on $g^{jk}$ is strong enough for
standard Schauder estimates to apply, yielding
$$
f_1\in C^r(\Ombar),\ {}\ h\in C^{r+2}(\Ombar)\Rightarrow u\in C^{r+2}(\Ombar),
\quad \text{for } r\in (0,1)\cup(1,2).
\tag{2.4.3}
$$
Then an interpolation argument gives the conclusion stated in Proposition 
2.4.1.
\enddemo

$\text{}$ \newline
{\smc Remark}. Under the hypotheses of Theorem 2.1, we also have
$$
g_{jk}\in H^{2,p}(\Omega),\quad \forall\ p<\infty,
\tag{2.4.4}
$$
in boundary harmonic coordinates.  Given this, we have, in parallel with
(2.4.1),
$$
\Delta u=f\in H^{1,q}(\Omega),\ u\bigr|_{\bW}=h\in B^{3-1/p}_{p,p}(\bW)
\Longrightarrow u\in H^{3,q}(\Omega),
\tag{2.4.5}
$$
for $q\in (1,\infty)$.  The proofs are simple variants of arguments presented 
above.

$\text{}$

\heading
\S{2.5}: Non-branching (and branching) of geodesics
\endheading

$\text{}$

Suppose $(g_{jk})$ is a metric tensor on an open set $\Cal{V}\subset\RR^n$.
The geodesic equation can be written in Hamiltonian form as
$$
\dot{x}^j=\frac{\pa}{\pa\xi_j}G(x,\xi),\quad 
\dot{\xi}_j=-\frac{\pa}{\pa x^j}G(x,\xi),
\tag{2.5.1}
$$
where
$$
G(x,\xi)=\frac{1}{2} g^{jk}(x)\xi_j \xi_k.
\tag{2.5.2}
$$
In other words, the geodesic flow is the flow on $T^*\Cal{V}=\Cal{V}\times
\RR^n$ generated by
$$
X(x,\xi)=\sum\limits_j \frac{\pa G}{\pa\xi_j}\frac{\pa}{\pa x^j}
-\frac{\pa G}{\pa x^j}\frac{\pa}{\pa \xi_j}.
\tag{2.5.3}
$$
Osgood's theorem states that a vector field generates a uniquely defined
flow provided its coefficients have a modulus of continuity $\omega(t)$
satisfying
$$
\int_0^{1/2} \frac{dt}{\omega(t)}=\infty.
\tag{2.5.4}
$$
See, e.g., Chapter 1 of [T1] for a proof.  An example for which (2.5.4)
holds is
$$
\omega(t)=t\, \log \frac{1}{t},
\tag{2.5.5}
$$
which is just a bit rougher than the Lipschitz modulus of continuity.
In fact, one has
$$
f\in C^1_*(U)\Longrightarrow |f(x+y)-f(x)|\le C \|f\|_{C^1_*}\, \omega(|y|),
\tag{2.5.6}
$$
with $\omega(t)$ given by (2.5.5).  See, e.g., [T2], Chapter I, for a proof
of this classical result.  This applies to (2.5.1)--(2.5.3) provided
$g_{jk}\in C^2_*(\Cal{V})$, so we have the following.

\proclaim{Proposition 2.5.1} Let $(g_{jk})\in C^2_*(\Cal{V})$ be a metric
tensor on an open set $\Cal{V}\subset\RR^n$.  Then the geodesic flow is
locally uniquely defined.
\endproclaim

\proclaim{Corollary 2.5.2} With $\Ombar$ as in Theorem 2.1, suppose
$g_{jk}\in C^2_*(\Ombar)$.  Then the geodesic flow is locally uniquely
defined when applied to any initial point $(z,\xi)$ with $z\in \Sigma$
and $v=v(\xi)\in T_z\Ombar$ pointing inside $\Omega$, transversal to 
$\Sigma$.
\endproclaim

Here $\xi$ and $v$ are related by $v^j=g^{jk}(z)\xi_k$.  To prove the 
corollary, let $\Cal{V}$ be a collar neighborhood of $\Ombar$.  One
can extend $g_{jk}$ to $g_{jk}\in C^2_*(\Cal{V})$, and we have that 
the geodesic flow is locally uniquely defined on $T^*\Cal{V}$.  If the
initial point is $(z,\xi)$, as described above, the geodesic has tangent
vector $v$ at $z$, so under our hypotheses the geodesic initially moves
into $\Omega$, for small positive $t$.  This behavior is hence independent
of the chosen extension of $g_{jk}$ to $\Cal{V}$.

We now discuss some examples of metric tensors only slightly rougher
than treated in Proposition 2.5.1, for which there is branching of geodesics.
These examples are variants of some produced by P.~Hartman in [Ha].
As in [Ha], we take
$$
ds^2=h(v)(du^2+dv^2).
\tag{2.5.7}
$$
As noted there, curves of the form $v=v(u)$ are (variable speed) geodesics 
provided
$$
2 \frac{d^2 v}{du^2}=\Bigl(1+\Bigl(\frac{dv}{du}\Bigr)^2\Bigr)H'(v),\quad
H(v)=\log h(v).
\tag{2.5.8}
$$
Multiplying by $dv/du$ yields
$$
\frac{d}{du}\Bigl(\frac{dv}{du}\Bigr)^2=\Bigl(1+\Bigl(\frac{dv}{du}\Bigr)^2
\Bigr) H'(v)\, \frac{dv}{du},
\tag{2.5.9}
$$
hence
$$
\frac{d}{du} \log (1+v^2_u)=\frac{d}{du} H(v(u)).
\tag{2.5.10}
$$
Thus $v=v(u)$ solves (2.5.8) provided $dv/du \neq 0$ and
$$
h(v(u))=1+\Bigl(\frac{dv}{du}\Bigr)^2.
\tag{2.5.11}
$$
If, however,
$$
h(0)=1,\quad h'(0)=0,
\tag{2.5.12}
$$
then $v\equiv 0$ solves (2.5.8), so if we produce another solution to (2.5.8) 
such that $v(0)=0,\ v'(0)=0$, we will have branching of geodesics.

For the first class of examples, we pick $k\in \{1,2,3,\dots\}$ and
construct $h(v)$ such that
$$
v(u)=u^{2k+1}
\tag{2.5.13}
$$
solves (2.5.8).  Indeed, by (2.5.11), this happens when
$$
h(v)=1+(2k+1)^2 |v|^{4k/(2k+1)}.
\tag{2.5.14}
$$
In this case the metric tensor has one derivative in $C^{(2k-1)/(2k+1)}$.
The case $k=1$ was explicitly mentioned in [Ha].

For another family of examples, we take $k\in\{1,2,3,\dots\}$ and 
construct $h(v)$ such that
$$
v(u)=e^{-1/|u|^k}\, \sgn u
\tag{2.5.15}
$$
solves (2.5.8).  This time we have
$$
u=\frac{\sgn v}{|\log|v||^{1/k}},\quad 
\Bigl(\frac{dv}{du}\Bigr)^2=
\frac{k^2}{u^{2k+2}}\, e^{-2/|u|^k},
\tag{2.5.16}
$$
and hence (2.5.11) gives
$$
h(v)=1+k^2 |\log|v||^{2+2/k}\, v^2.
\tag{2.5.17}
$$
For these examples, the first derivative of the metric tensor has a 
modulus of continuity only slightly worse than log-Lipschitz.
Also the metric tensor has two derivatives in $L^p$ for all $p<\infty$.
Hence the curvature belongs to $L^p$ for all $p<\infty$.

$$\text{}$$
{\bf 3. Geometric convergence for manifolds with boundary}
\newline {}\newline

A sequence $(\Mbar_k,g_k)$ of compact Riemannian manifolds with boundary
$\pa M_k$ is said to converge in the $C^r$-topology (given $0<r<\infty$) to a
compact Riemannian manifold $(\Mbar,g)$ provided that $g$ is a $C^r$ metric
tensor on $\Mbar$ and, for $k$ sufficiently large, there exist 
diffeomorphisms $F_k:\Mbar\rightarrow \Mbar_k$ such that $F^*_kg_k$ 
converges to $g$ in the $C^r$-topology.  (Necessarily $F_k:\pa M\rightarrow
\pa M_k$.)  In this section we will identify classes of Riemannian manifolds 
with boundary that are pre-compact in the $C^r$-topology, for any given
$r<2$.

We work with families of Riemannian manifolds with boundary of the following
sort.  Fix the dimension, $n$.
Given $R_0,i_0,S_0,d_0\in (0,\infty)$, denote by $\Cal{M}(R_0,i_0,
S_0,d_0)$ the class of compact, connected, $n$-dimensional
Riemannian manifolds with boundary $(\Mbar,g)$, with smooth metric tensor,
with the following four properties:
$$
\|\Ric_M\|_{L^\infty(M)}\le R_0,\quad \|\Ric_{\pa M}\|_{L^\infty(\pa M)}
\le R_0, 
\tag{3.0.1}
$$
where $\Ric$ denotes the Ricci tensor.
$$
i_M\ge i_0,\quad i_{\pa M}\ge i_0,\quad i_b\ge 2i_0.
\tag{3.0.2}
$$
Here $i_M$ denotes the injectivity radius of $\Mbar$, $i_{\pa M}$ that of
$\pa M$, and $i_b$ the boundary injectivity radius of $\Mbar$.
$$
\|H\|_{\Lip(\pa M)}\le S_0,
\tag{3.0.3}
$$
where $H$ is the mean curvature of $\pa M$ in $\Mbar$.
$$
\text{diam}\, (\Mbar,g)\le d_0.
\tag{3.0.4}
$$

We recall the concept of boundary injectivity radius, $i_b$.
It is the optimal quantity with the following property. 
Namely, there is a collar neighborhood
$\Cal{C}$ of $\pa M$ in $\Mbar$ and a (unique) function $f\in C^2(\Cal{C})$
such that
$f\bigr|_{\pa M}=0,\quad |\nabla f|\equiv 1,\quad
f(\Cal{C})\supset [0,i_b).$
With this, local coordinates $(v^1,\dots,v^{n-1})$ on an open set in $\pa M$
can be continued inside, as constant on the integral curves of $\nabla f$, to 
produce, along with $v^n=f$, a set of ``boundary normal coordinates.''

To further clarify the first part of (3.0.2), we mean that
$$
\Exp_p:B_\rho(0)\rightarrow M,
\tag{3.0.5}
$$
where $B_\rho(0)=\{v\in T_pM:g(v,v)<\rho^2\}$, is a diffeomorphism for
$\rho=i_0$ if $\dist(p,\pa M)\ge i_0$ and it is a diffeomorphism for $\rho=
\dist(p,\pa M)$ if $\dist(p,\pa M)\le i_0$.

The main goal in this section is to prove the following.

\proclaim{Theorem 3.1} Given $R_0,i_0,S_0,d_0\in (0,\infty)$, $\Cal{M}
(R_0,i_0,S_0,d_0)$ is precompact in the $C^r$-topology for each $r<2$.  
In particular, any sequence $(\Mbar_k,g_k)$ in $\Cal{M}(R_0,i_0,S_0,d_0)$
has a subsequence that converges in the $C^r$-topology to a limit 
$(\Mbar,g)$.  Furthermore, the metric tensor $g$ belongs to $C^2_*(\Mbar)$.
\endproclaim

Such a result was established in [An1] in the category of compact manifolds
without boundary; subsequently there have been expositions in [HH] and in
[Pe].  Our proof of Theorem 3.1 follows the structure of the argument in
[An1], with necessary modifications to treat the case of nonempty 
boundary.  In this regard the boundary regularity results of \S{2.1} play
a major role.  The $C^2_*$ part of the conclusion is also more precise than
that noted in earlier results.  This precision will be of major value in
the application of Theorem 3.1 to results on inverse boundary spectral 
problems presented in \S{4}.

The proof of Theorem 3.1 involves a blow-up argument that takes us outside
the category of compact manifolds, and it is useful to have the following
notion of pointed convergence of a sequence $(\Mbar_k,g_k,p_k)$, with 
distinguished points $p_k\in\Mbar_k$.  We say $(\Mbar_k,g_k,p_k)$ converges
to $(\Mbar,g,p)$ in the pointed $C^r$-topology provided the following holds.
For large $k$ there exist $\rho_k < \sigma_k, \rho_k
\nearrow  \infty $, 
and compact 
$\Ombar_k\subset \Mbar_k$ and $\overline{\Cal{V}}_k\subset\Mbar$, such that
$$
B_{\rho_k}(p_k)\subset \Ombar_k\subset B_{\sigma_k}(p_k),\quad 
B_{\rho_k}(p)\subset \overline{\Cal{V}}_k\subset B_{\sigma_k}(p),
$$
and diffeomorphisms 
$$
F_k:\overline{\Cal{V}}_k\rightarrow \Ombar_k,\quad
F_k:\overline{\Cal{V}}_k\cap \pa M\rightarrow \Ombar_k\cap \pa M_k,
$$
such that
$F_k^*g_k$ converges to $g$ in the $C^r$-topology (on each compact
subset of $\Mbar$) and $F_k^{-1}(p_k)\rightarrow p$.  We assume $\Mbar$
is connected, and these hypotheses imply $(\Mbar,g)$ must be complete.

$\text{}$

\heading
\S{3.1}: Basic convergence results
\endheading

$\text{}$

We begin by describing a slight variant of a well known ``abstract'' 
convergence result.  As before, fix the dimension $n$.  Given $s=\ell+\sigma$
($\ell\in\ZZ^+,\ 0<\sigma<1$), $\rho\in(0,\infty),\ Q\in (1,2)$, let $\Cal{N}
(s,\rho,Q)$ denote the class of connected, $n$-dimensional Riemannian 
manifolds $(\Mbar,g)$ with boundary, with the following properties.
Take $p\in\Mbar$.  
\newline {}\newline
(i) If $\text{dist}\,(p,\pa M)>\rho$, there is a neighborhood $U$ of $p$
in $M^{int}$ and a coordinate chart $\varphi:B_{\rho/2}(0)\rightarrow U$, such
that, in these coordinates
$$
Q^{-2}|\eta|^2\le g_{jk}(x)\eta^j\eta^k\le Q^2|\eta|^2,
\tag{3.1.1}
$$
and
$$
\rho^{s} \sum\limits_{|\beta|=\ell} \sup |x-y|^{-\sigma}\, 
|\pa^\beta g_{jk}(x)-\pa^\beta g_{jk}(y)|\le Q-1.
\tag{3.1.2}
$$
(ii) If $\text{dist}\,(p,\pa M)\le \rho$, there is a neighborhood $U$ of $p$
in $\Mbar$ and a coordinate chart $\varphi:B^+_{4\rho}(0)\rightarrow U$ 
such that $\{x^n=0\}$ maps to $\pa M$ and (3.1.1)--(3.1.2) hold in these
coordinates.
\newline {}\newline
Let $\Cal{N}_*(s,\rho,Q)$ denote the class of pointed manifolds
$(\Mbar,g,p)$, with $p\in\Mbar$, satisfying these properties.
\newline $\text{}$

The following compactness result goes back to [Ch]; a detailed proof is given
on pp.~293--296 of [Pe] for the case of manifolds without boundary, but no
essential changes are required for the case of manifolds with boundary.

\proclaim{Theorem 3.1.1} Given $s,\rho\in (0,\infty),\ Q\in (1,2)$, the class 
$\Cal{N}_*(s,\rho,Q)$ is compact in $\Cal{N}_*(s',\rho,Q)$ in the pointed
$C^{s'}$-topology, for all $s'<s$.
\endproclaim

As for convergence without specifying base points, we can define 
$\Cal{N}(s,\rho,Q,d_0)$ to consist of $(\Mbar,g)\in \Cal{N}(s,\rho,Q)$ 
satisfying
also $\text{diam}\,(\Mbar,g)\le d_0$, and conclude that $\Cal{N}(s,\rho,Q,d_0)$
is compact in $\Cal{N}(s',\rho,Q,d_0)$ for all $s'<s$.

To apply Theorem 3.1.1 to our situation, we will show that, given $R_0,i_0,
S_0,d_0\in (0,\infty)$, and given $Q\in (1,2),\ s\in (1,2)$, there exists
$\rho>0$ such that
$$
\Cal{M}(R_0,i_0,S_0,d_0)\subset \Cal{N}(s,\rho,Q).
\tag{3.1.3}
$$
In fact, we will establish a result that is more precise, in two respects.

For one, we will produce harmonic coordinates, in case (i), and boundary 
harmonic coordinates, in case (ii).  Doing this brings in the notion of
$C^s$-harmonic radius, introduced in [An1] in the context of manifolds
without boundary.  To be precise, if $s=\ell+\sigma$, as above, a number
$r^s_h=r^s_h(p,g,Q)$ is called the $C^s$-harmonic radius of $(\Mbar,g)$
at $p\in\Mbar$ provided it is the optimal quantity with the following 
property.  Take any $\rho<r^s_h$.  Then, in case (i), there exist 
harmonic coordinates $\varphi^{-1}:U\rightarrow B_{\rho/2}(0)$ such that
(3.1.1)--(3.1.2) hold on $B_{\rho/2}(0)$.
In case (ii), there exist harmonic coordinates
$\varphi^{-1}:U\rightarrow B^+_{4\rho}(0)$, such that (3.1.1)--(3.1.2) hold
on $B^+_{4\rho}(0)$.  If $\Mbar$ is compact, we define the $C^s$-harmonic
radius of $(\Mbar,g)$ by
$$
r^s_h(\Mbar,g,Q)=\inf\limits_{p\in\Mbar} \, r^s_h(p,g,Q).
\tag{3.1.4}
$$
The containment (3.1.3) follows from the fact that, given $R_0,i_0,S_0,d_0
\in (0,\infty)$ and $s\in (1,2),\, Q \in (1,2)$, there is a lower bound on the 
$C^s$-harmonic radius of $(\Mbar,g)\in\Cal{M}(R_0,i_0,S_0,d_0)$.  

The result we will establish is more precise in one more respect;
we bring in the notion of $C^2_*$-harmonic radius.  Namely, a number
$r_h=r_h(p,g,Q)$ is called the $C^2_*$-harmonic radius of $(\Mbar,g)$
at $p\in\Mbar$ given the circumstances described above (in the definition
of $C^s$-harmonic radius) but with (3.1.2) replaced by
$$
\rho^2\sum\limits_{|\beta|=1} \sup\, |x-y|^{-1} \bigl|\pa^\beta g_{jk}(x)
+\pa^\beta g_{jk}(y)-2\pa^\beta g_{jk}\bigl((x+y)/2\bigr)\bigr|\le Q-1.
\tag{3.1.5}
$$
Then the $C^2_*$-harmonic radius of $(\Mbar,g)$ is defined by
$$
r_h(\Mbar,g,Q)=\inf\limits_{p\in\Mbar}\, r_h(p,g,Q).
\tag{3.1.6}
$$

We note the following elementary but important scaling property of 
harmonic radius, valid for $\lambda\in (0,\infty)$:
$$
r_h(\Mbar,\lambda^2 g,Q)=\lambda r_h(\Mbar,g,Q),\quad
r^s_h(\Mbar,\lambda^2 g,Q)=\lambda r^s_h(\Mbar,g,Q).
\tag{3.1.7}
$$

Our next goal is to show that, given $R_0,i_0,S_0,d_0\in (0,\infty)$, 
there is a lower bound on the $C^2_*$-harmonic radius of $(\Mbar,g)\in
\Cal{M}(R_0,i_0,S_0,d_0)$.

$\text{}$

\heading
\S{3.2}: Harmonic radius estimate
\endheading

$\text{}$

As advertised, our goal here is to prove the following.

\proclaim{Theorem 3.2.1} Let $R_0, i_0,S_0$, and $d_0$ be given, in 
$(0,\infty)$, and let $Q\in (1,2)$ be given.  
Then there exists $r_{\Cal{M}}=r_{\Cal{M}}(R_0,i_0,
S_0,d_0,Q)>0$ such that
$$
r_h(\Mbar,g,Q)\ge r_{\Cal{M}},\quad \forall\ (\Mbar,g)\in 
\Cal{M}(R_0,i_0,S_0,d_0).
\tag{3.2.1}
$$
\endproclaim

The proof will be by contradiction.  Suppose there exist 
$(\Mbar_k,\tilde{g}_k)\in\Cal{M}(R_0,i_0,S_0,d_0)$ such that
$$
r_h(\Mbar_k,\tilde{g}_k,Q)=\ep_k\rightarrow 0.
\tag{3.2.2}
$$
Let us scale the metric $\tilde{g}_k$ to $g_k=\ep^{-2}_k\tilde{g}_k$, 
and consider the scaled Riemannian manifolds $(\Mbar_k,g_k)$.  Then
$$
r_h(\Mbar_k,g_k,Q)=1,
\tag{3.2.3}
$$
while, for the rescaled metric,
$$
\|\Ric_k\|_{L^\infty(M)}\rightarrow 0,\quad \|\Ric_{\pa M_k}
\|_{L^\infty(\pa M_k)}\rightarrow 0,\quad \|H_k\|_{\Lip^1(\pa M_k)}
\rightarrow 0,
\tag{3.2.4}
$$
and
$$
i_k\rightarrow\infty,\quad i_{\pa M_k}\rightarrow\infty,\quad i_{b,M_k}
\rightarrow \infty,
\tag{3.2.5}
$$
as $k\rightarrow\infty$.  We will show that conditions (3.2.3)--(3.2.5) 
lead to a contradiction.

To see this, pick $p_k\in M_k$ such that 
$$
r_h(p_k,g_k,Q)=1.
\tag{3.2.6}
$$
By Theorem
3.1.1, there is a subsequence, which for simplicity we will also denote 
$(\Mbar_k,g_k,p_k)$ which converges to a pointed manifold $(\Mbar,g,p)$
in the $C^r$-topology, where $r\in (1,2)$ is arbitrary.  
With $\tau_k = \text{dist}_k(x, \partial M_k)$, where
$\text{dist}_k$ is the distance on $(\Mbar_k,g_k)$,
there are two possibilities:
\newline {}\newline
(i) $\tau_k(p_k)\rightarrow\infty$.

Then we claim $(M,g)$ is isometric to $\RR^n$, with its standard flat metric.
\newline {}\newline
(ii) $\tau_k(p_k)\le K<\infty$.

Then we claim $(\Mbar,g)$ is isometric to $\overline{\RR}^n_+$, with its
standard flat metric.  
\newline $\text{}$

We demonstrate these claims in Lemma 3.2.2 below.  For now we assume this, 
and proceed with a further analysis of these two cases.

$\text{}$ \newline
{\smc Case (i)}.  This case is treated in [An1].  We recall the argument
here, in a slightly varied form (also borrowing from [HH]),
both to establish a slightly stronger conclusion and to set the stage
to examine Case (ii).

We have neighborhoods $U_k$ of $p_k$ in $M_k$, identified with $B_5=
\{x\in\RR^n:|x|\le 5\}$, with $p_k=0$, and the metric tensors $g_k\rightarrow
\delta$ in $C^r$-norm on $B_5$, where $\delta$ is the standard Euclidean
metric tensor on $\RR^n$.  Taking $x^\nu\ (1\le\nu\le n)$ to be the standard
Cartesian coordinates on $\RR^n$, we solve
$$
\Delta_k u^\nu_k=0\ \text{ on }\ B_5,\quad u^\nu_k=x^\nu\ \text{ on }\
\pa B_5,
\tag{3.2.7}
$$
where $\Delta_k$ is the Laplace operator with respect to the metric tensor
$g_k$:
$$
\Delta_k u=g_{(k)}^{-1/2} \pa_i \bigl(g^{1/2}_{(k)}g_{(k)}^{ij} \pa_j u
\bigr).
\tag{3.2.8}
$$
Note that
$$
\Delta_k (u^\nu_k-x^\nu)=f^\nu_k\ \text{ on }\ B_5,\quad
u^\nu_k-x^\nu=0\ \text{ on }\ \pa B_5,
\tag{3.2.9}
$$
where
$$
f^\nu_k=-g^{-1/2}_{(k)} \pa_i(g^{1/2}_{(k)} g^{i\nu}_{(k)}),
\tag{3.2.10}
$$
for which we have
$$
\|f^\nu_k\|_{C^{r-1}(B_5)}\rightarrow 0,\quad k\rightarrow\infty.
\tag{3.2.11}
$$
In view of this and the uniform $C^r$ estimates on the coefficients in the 
elliptic differential operators (3.2.8), elliptic regularity gives
$$
u^\nu_k-x^\nu\rightarrow 0\ \text{ in }\ C^{r+1}(B_5),\quad
k\rightarrow\infty.
\tag{3.2.12}
$$
Hence, for all large $k$, $(u^1_k,\dots,u^n_k)$ form a coordinate 
system on $U_k$, and if now $g^{(k)}_{ij}$ denotes the components of the
metric tensor $g_k$ in this coordinate system, then
$$
\|g^{(k)}_{ij}-\delta_{ij}\|_{C^r(B_4)}\rightarrow 0,\quad 
k\rightarrow\infty.
\tag{3.2.13}
$$

The last step is to obtain an analogue of (3.2.13) in a stronger norm, using 
the Ricci equation, which implies
$$
\Delta_k (g^{(k)}_{ij}-\delta_{ij})=F^{(k)}_{ij},
\tag{3.2.14}
$$
where $\Delta_k$ is as in (3.2.8), acting componentwise, and
$$
F^{(k)}_{ij}=B_{ij}(g_{\ell m},\nabla g_{\ell m})-2\Ric^{(k)}_{ij}.
\tag{3.2.15}
$$
Here $B_{ij}$ is a quadratic form in $\nabla g_{\ell m}$.  In view of (3.2.4)
and (3.2.13) we have $\|F^{(k)}_{ij}\|_{L^\infty(B_4)}\rightarrow 0$.  Hence 
local elliptic regularity results on (3.2.14) (plus another appeal to 
(3.2.13)) give
$$
\|g^{(k)}_{ij}-\delta_{ij}\|_{C^2_*(B_3)}\rightarrow 0,\quad 
k\rightarrow\infty,
\tag{3.2.16}
$$
and even the stronger result
$$
\|g^{(k)}_{ij}-\delta_{ij}\|_{\frak{h}^{2,\infty}(B_3)}\rightarrow 0,\quad
k\rightarrow\infty.
\tag{3.2.17}
$$

Hence we have $r_h(p_k,g_k,Q)\ge 3$ for large $k$, contradicting (3.2.6).
\newline {}\newline
{\smc Remark}. Clearly there is nothing special about ``3'' in (3.2.16).
For any $L\in (1,\infty)$ one obtains (3.2.16) with $C^2_*(B_3)$ replaced 
by $C^2_*(B_L)$.

$\text{}$ \newline
{\smc Case (ii)}.  In this case, for large $k$, there is a unique 
$q_k\in\pa M_k$ closest to $p_k$.  In addition to the pointed convergence
$(\Mbar_k,g_k,p_k)\rightarrow (\Mbar,g,p)$, we also have $(\pa M_k,g_k,q_k)
\rightarrow (\pa M,g,q)$, and $\text{dist}_k(p_k,q_k)\rightarrow 
\text{dist}(p,q)$.  
This time $(\Mbar,\pa M,g)=(\overline{\RR}^n_+,\RR^{n-1},\delta)$, with
$\RR^{n-1}=\{x^n=0\}$, and we can identify $q$ with $0\in\RR^{n-1}\subset
\RR^n$.  Fix $L\ge 2\, \text{dist}(p,q)+4$.

In this case we have neighborhoods $U_k$ of $q_k$ in $\Mbar_k$, 
identified with
$$
B^+_{L+5}=\{x\in\overline{\RR}^n_+:|x|\le L+5\},
$$
with $q_k=0$, and the metric tensors $g_k\rightarrow \delta$ in $C^r$-norm
on $B^+_{L+5}$.  Note that 
$\{x\in M_k:\text{dist}_k(x,p_k)\le 2\}\subset B^+_{L+5}$,
for large $k$.

With $x^\nu$ ($1\le\nu\le n$) the standard Cartesian coordinates on 
$\overline{\RR}^n_+$ ($x^n=0$ defining $\pa \overline{\RR}^n_+$) we first 
solve, for $1\le\nu\le n-1$,
$$
\Delta_{\pa M_k} v^\nu_k=0\ \text{ on }\ \widetilde{B}_{L+5},\quad
v^\nu_k=x^\nu\ \text{ on }\ \pa\widetilde{B}_{L+5},
\tag{3.2.18}
$$
where
$$
\widetilde{B}_{L+5}=\{x\in\RR^{n-1}:|x|\le L+5\}.
$$
Parallel to (3.2.12), we have
$$
\|v^\nu_k-x^\nu\|_{C^{r+1}(\widetilde{B}_{L+5})}\rightarrow 0,\quad
k\rightarrow 0,\quad 1\le\nu\le n-1.
\tag{3.2.19}
$$
Next solve
$$
\Delta_k u^\nu_k=0\ \text{ on }\ B^+_{L+5},\quad
u^\nu_k\bigr|_{\widetilde{B}_{L+5}}=v^\nu_k,\quad
u^\nu_k\bigr|_{\pa^+B_{L+5}}=x^\nu,
\tag{3.2.20}
$$
where
$$
\pa^+ B_{L+5}=\pa(B^+_{L+5})\setminus \widetilde{B}_{L+5}.
$$
Parallel to (3.2.9)--(3.2.11), we have
$$
\Delta_k(u^\nu_k-x^\nu)=f^\nu_k\ \text{ on }\ B^+_{L+5},\quad
u^\nu_k-x^\nu=\varphi^\nu_k\ \text{ on }\ \pa B^+_{L+5},
\tag{3.2.21}
$$
with
$$
\|f^\nu_k\|_{C^{r-1}(B^+_{L+5})}\rightarrow 0,
\tag{3.2.22}
$$
and
$$
\|\varphi^\nu_k \|_{L^\infty(\pa B^+_{L+5})}\rightarrow 0,
\quad \|\varphi^\nu_k\|_{C^{r+1}(\widetilde{B}_{L+5})}\rightarrow 0.
\tag{3.2.23}
$$
This gives first a global estimate
$$
\|u^\nu_k-x^\nu\|_{L^\infty(B^+_{L+5})}\rightarrow 0,
\tag{3.2.24}
$$
and then a ``local'' estimate (away from the corners of $B^+_{L+5}$)
$$
\|u^\nu_k-x^\nu\|_{C^{r+1}(B^+_{L+4})}\rightarrow 0,
\tag{3.2.24}
$$
as $k\rightarrow\infty$, for $1\le\nu\le n-1$.

To construct $u^n_k$, we solve
$$
\Delta_k u^n_k=0\ \text{ on }\ B^+_{L+5},\quad
u^n_k=x^n\ \text{ on }\ \pa B^+_{L+5}.
\tag{3.2.26}
$$
Parallel to (3.2.25) we have
$$
\|u^n_k-x^n\|_{C^{r+1}(B^+_{L+4})}\rightarrow 0,\quad k\rightarrow\infty.
\tag{3.2.27}
$$
Hence, for large $k,\ (u^1_k,\dots,u^n_k)$ form a coordinate system
on $B^+_{L+4}$, with $\{u^n_k=0\}$ defining the face $\{x^n=0\}$.
Parallel to (3.2.13), if $g^{(k)}_{ij}$ denotes the components of the metric
tensor $g_k$ in this coordinate system, we have
$$
\|g^{(k)}_{ij}-\delta_{ij}\|_{C^r(\pa B^+_{L+3})}\rightarrow 0,
\quad k\rightarrow\infty,
\tag{3.2.28}
$$
for any $r<2,\ 1\le i,j\le n$.

As in Case (i), the last step is to obtain an analogue of (3.2.28) in a 
stronger norm.  First, if we set $h^{(k)}_{ij}=g^{(k)}_{ij}
\bigr|_{\widetilde{B}_{L+5}},
\,1\le i,j\le n-1$, then the argument in Case (i) applies to yield
$$
\|h^{(k)}_{ij}-\delta_{ij}\|_{\frak{h}^{2,\infty}(\widetilde{B}_{L+2})}
\rightarrow 0,\quad k\rightarrow\infty.
\tag{3.2.29}
$$
Now estimates such as derived in \S{2} apply to $g^{(k)}_{ij}$, solving
(3.2.14)--(3.2.15) and the boundary condition (3.2.29), for 
$1\le i,j\le n-1$, to yield
$$
\|g^{(k)}_{ij}-\delta_{ij}\|_{C^2_*(B^+_{L+1})}\rightarrow 0,
\quad k\rightarrow \infty,
\tag{3.2.30}
$$
for $1\le i,j\le n-1$.  Next, we use
$$
\Delta_k(g^{\ell n}_{(k)}-\delta^{\ell n})=F^{\ell n}_{(k)}
=B^{\ell n}(g_{(k)},\nabla g_{(k)})+2 \, \Ric_{(k)}^{\ell n},
\tag{3.2.31}
$$
with boundary conditions
$$
\aligned
N(g^{nn}_{(k)}-\delta^{nn})&=-2(n-1)H_kg^{nn}_{(k)}, \\
N(g^{\ell n}_{(k)}-\delta^{\ell n})&=-(n-1)H_k g^{\ell n}_{(k)}
+\frac{1}{2} \frac{1}{\sqrt{g^{nn}_{(k)}}} g^{\ell i}_{(k)}
\pa_i g^{nn}_{(k)},\quad 1\le\ell\le n-1,
\endaligned
\tag{3.2.32}
$$
on $\widetilde{B}_{L+5}$.  The results of \S{2} together with hypotheses 
in (3.2.4) and the analogue of (3.2.28) for $\|g^{ij}_{(k)}-\delta^{ij}
\|_{C^r(B^+_{L+3})}$, yield first
$$
\|g^{nn}_{(k)}-\delta^{nn}\|_{C^2_*(B^+_L)}\rightarrow 0,\quad
k\rightarrow\infty,
\tag{3.2.33}
$$
and then
$$
\|g^{\ell n}_{(k)}-\delta^{\ell n}\|_{C^2_*(B^+_L)}\rightarrow 0,
\quad k\rightarrow\infty,\quad \ell<n.
\tag{3.2.34}
$$
Finally an argument parallel to that in Lemma 2.1.3 gives (3.2.30) for all
$i,j\le n$.

Hence we have $r_h(p_k,g_k,Q) \geq 2$ for large $k$, contradicting (3.2.6).
This finishes the proof of Theorem 3.2.1, modulo the proof of the following
lemma.

\proclaim{Lemma 3.2.2} Assume that a sequence of Riemannian manifolds
$(\Mbar_k,g_k,p_k)$ satisfies the conditions of Theorem 3.2.1 and 
$(\Mbar_k,g_k,p_k)\rightarrow (\Mbar,g,p)$ 
in the $C^{r}$-topology, $r \in(1,2)$.  Assume that
$$
\|\Ric_k\|_{L^\infty(M_k)}\rightarrow 0,\quad
\|\Ric_{\pa M_k}\|_{L^\infty(\pa M_k)}\rightarrow 0,\quad
\|H_k\|_{L^\infty(\pa M_k)}\rightarrow 0,
\tag{3.2.35}
$$
and
$$
i_{M_k}\rightarrow\infty,\quad i_{\pa M_k}\rightarrow \infty,\quad
i_{b,M_k}\rightarrow\infty,
\tag{3.2.36}
$$
as $k\rightarrow\infty$.  Then $(\Mbar,g)$ is isometric to either $\RR^n$
or $\overline{\RR}^n_+$.
\endproclaim
\demo{Proof} We consider 
separately the cases $\tau_k(p_k)\rightarrow\infty$ and $\tau_k(p_k)\le K
<\infty$.  
\newline {}\newline
{\smc Case (i)}. Suppose $\tau_k(p_k)\rightarrow\infty$. 

This case is effectively treated in [An1].  We recall briefly the argument,
since it also plays a role in Case (ii).  In this case the hypotheses imply
$$
\Ric_M=0,
\tag{3.2.37}
$$
weakly.  Hence the metric tensor $g$ is smooth in local harmonic coordinates.
Also any unit speed geodesic $\gamma(t)$ such that $\gamma(0)=p$
is defined in $M$ for all $t\in\RR$.  Take $T\in (0,\infty)$.  We claim 
$\gamma$ is the shortest path from $\gamma(-T)$ to $\gamma(T)$.

Consider $B_{6T}(p)=\{x\in M:\text{dist}(x,p)<6T\}$.
The hypotheses imply that there exists $k_0$ such that for all $k\ge k_0$
there are open sets in $M_k$ identified with $B_{6T}(p)$ via diffeomorphisms,
such that $p_k\rightarrow  p$ and $g_k\rightarrow g$ in $C^{r}(B_{6T}(p))$.
Also, for each $x\in B_{2T}(p),\ i(x,g_k)\ge 3T$.

Since $(M,g)$ is smooth, there is some $c_0>0$ such that $i(p,g)\ge c_0$.
Say $\gamma(c_0)=q$.  We have unit speed geodesics $\gamma_k$ on $(M_k,g_k)$,
defined for $|t|\le 2T$, such that $\gamma_k(0)=p,\ \gamma_k(c_k)=q,\
c_k\rightarrow c_0$.  By Arzela's theorem there is a subsequence $\gamma_k
\rightarrow \sigma$, uniformly on $t\in[-2T,2T]$.  We see that $\sigma$ is a 
geodesic on $(M,g),\ \sigma(0)=p,\ \sigma(c_0)=q$.  Hence $\sigma=\gamma$.
(It follows that the entire sequence $\gamma_k\rightarrow\gamma$, not just a 
subsequence.)  Since $i(x,g_k)\ge 3T$ for $x\in B_{2T}(p)$, we can deduce
that
$$
\text{dist}_k(\gamma_k(-T),\gamma_k(T))=2T,
\tag{3.2.38}
$$
for $k\ge k_0$.  Since $\gamma_k\rightarrow \gamma$ uniformly on $[-T,T]$
and $g_k\rightarrow g$ in $C^0$, the left side of (3.2.38) converges to 
$\text{dist}(\gamma(-T),\gamma(T))$.

Thus $(M,g)$ is complete and Ricci flat and each geodesic through $p$ is
globally length minimizing.  It follows from the Cheeger-Gromoll splitting
theorem [CG] that $(M,g)$ is isometric to standard flat $\RR^n$.

$\text{}$\newline
{\smc Case (ii)}.  Assume $\tau_k(p_k)\le K<\infty$.
\newline {}\newline

Let $q_k\in\pa M_k$ be the nearest point on $\pa M_k$ to $p_k$.  By (3.2.36)
$q_k$ is uniquely defined (for large $k$).  Also $(\pa M_k,g_k,q_k)$
converges in the $C^{r}$ topology to $(\pa M,g,q)$, and
$\text{dist}(p,q)=\lim \tau_k(p_k)$.  In this case the hypotheses yield
$$
\Ric_M=0,\quad \Ric_{\pa M}=0,\quad H\bigr|_{\pa M}=0.
\tag{3.2.39}
$$
Hence $g$ is smooth on $\Mbar$ in local harmonic coordinates.  

We will make strong use of two equations.  One is:
$$
\Delta_k \tau_k=H_k\ \text{ on }\ \{\tau_k=c\}.
\tag{3.2.40}
$$
Here, in slight contrast to notation in \S{2}, we set $H_k$ equal to
the trace of the Weingarten map, i.e., to $n-1$ times
the mean curvature of the surface $\{\tau_k=c\}$.  The other
is (in boundary normal coordinates $(z,\tau_k),\, z \in \pa M$): 
$$
\pa_\tau H_k=-\Tr A^2_k-\Ric_{(k)}(\pa_\tau,\pa_\tau),
\tag{3.2.41}
$$
where $A_k$ is the Weingarten map of the surface $\{\tau_k=c\}$.
See [Pe], \S\S{2.3--2.4}.

Under our hypotheses, there exist $\ep_k\rightarrow 0$ such that
$$
\|\Ric_{(k)}(\pa_\tau,\pa_\tau)\|_{L^\infty}\le \ep_k,\quad 
\|H_k(0)\|_{L^\infty}\le\ep_k.
\tag{3.2.42}
$$
Hence
$$
\pa_\tau H_k\le -\frac{1}{n-1} H^2_k+\ep_k.
\tag{3.2.43}
$$
In particular, $H_k(\tau,z)\le \mu_k(\tau)$ for $0\le\tau\le i_{b,M_k}$,
where
$$
\pa_\tau \mu_k=-\frac{1}{n-1}\mu_k^2+\ep_k,\quad \mu_k(0)=\ep_k.
\tag{3.2.44}
$$
Hence
$$
H_k(\tau,z)\le\max(\sqrt{(n-1)\ep_k},\ep_k),
\tag{3.2.45}
$$
since (3.2.44) forces $\mu_k'(\tau)\le 0$ before $\mu_k(\tau)$
can be larger than the right side of (3.2.45).

We next estimate $H_k(\tau,z)$ from below.  Pick $\delta>0$ and suppose
$$
H_k(\tau_0,z)\le-\delta,\ \text{for some }\tau_0\le\frac{i_{b,M_k}}{2},\
z\in\pa M.
\tag{3.2.46}
$$
For short, we set $H_k(\tau)=H_k(\tau,z)$.  If $k$
is large enough that $\delta^2>2(n-1)\ep_k$, it follows that $H_k(\tau)<
-\delta$ for all $\tau\in[\tau_0,i_{b,M_k})$, and hence
$$
\pa_\tau H_k\le -\frac{1}{2(n-1)}H^2_k,\quad \tau\ge\tau_0.
\tag{3.2.47}
$$
As $H_k(\tau_0)\le -\delta$, this implies
$$
H_k(\tau)\le \frac{2(n-1)}{(\tau-\tau_0)-2(n-1)/\delta}.
\tag{3.2.48}
$$
But this implies blow-up of $H_k(\tau)$ somewhere on $\tau\in[\tau_0,\tau_0
+2(n-1)/\delta]$, contradicting the fact that $H_k(\tau)$ is finite on
$\tau\in[0,i_{b,M_k})$, which contains $[\tau_0,\tau_0+2(n-1)/\delta]$
when $\delta >4(n-1)/i_{b,M_k}$. 
This contradiction shows that, for any given $\delta>0$, (3.2.46) must fail
for all sufficiently large $k$.  Hence we have $\delta_k\rightarrow 0$ 
such that
$$
|H_k(\tau)|\le\delta_k\ \text{ for }\ 0\le\tau\le \frac{i_{b,M_k}}{2}.
\tag{3.2.49}
$$
If $\tau(x)=\text{dist}(x,\pa M)$, then $\tau_k \rightarrow \tau$
in $C^0$ due to $g_k\rightarrow g$ in $C^{r}$.
Using again that $g_k\rightarrow g$ in $C^{r}$ together with
(3.2.40), (3.2.49), the boundary condition  
$ \tau_k|_{\pa M_k}=0$, and Propositions 5.1.1 and
5.2.1, we deduce that, for all $ s<2, \, n/(n-1) \le p < \infty$,
$$
\tau_k\rightarrow \tau\ \text{ in }\ C^{s}\cap H^{2,p}(B^+(\rho_k)), 
\tag{3.2.50}
$$
where $\rho_k\rightarrow\infty$ for $k\rightarrow\infty$.  In particular,
$i_{b,M}=+\infty$. Applying again (3.2.40), (3.2.49) together with
$g_k\rightarrow g$ in $C^{r}$, we obtain that
$$
H_k \rightarrow H \quad \text{in} \, \, L^p(B^+(\rho_k)),
\quad  \frac{n}{n-1} \leq p < \infty,
\tag{3.2.51}
$$
i.e., 
$$
H(\tau)=0,\quad \forall\ \tau\ge 0.
\tag{3.2.52}
$$

Now, parallel to (3.2.41), we have
$$
\pa_\tau H+\Tr A^2=-\Ric_M(\pa_\tau,\pa_\tau)=0,
\tag{3.2.53}
$$
where $A$ is the Weingarten map on the surface $\{\tau=c\}$,
and hence $\Tr A^2=0$, which implies
$$
A=0.
\tag{3.2.54}
$$
Furthermore, since
$$
\pa_\tau g_{ij}=2A^\ell{}_i\, g_{\ell j}=0,
\tag{3.2.55}
$$
we have
$$
g_{ij}(\tau,z)=g_{ij}(0,z),\quad z\in\pa M.
\tag{3.2.56}
$$
By the argument of Case (i),
$\pa M$ is isometric to $\RR^{n-1}$, so this implies that $\Mbar$
is isometric to $[0,\infty)\times\pa M=\overline{\RR}^n_+$, with its
standard flat metric.
\enddemo

This finishes the proof of Lemma 3.2.2, hence of Theorem 3.2.1.

$\text{}$

\heading
\S{3.3}: Proof of Theorem 3.1
\endheading

$\text{}$

All the work needed to prove Theorem 3.1 has been done, and we need only 
collect the pieces.  Say $(\Mbar_k,g_k)\in\Cal{M}(R_0,i_0,S_0,d_0)$.
The results of \S{3.2} show that Theorem 3.1.1 is applicable.  Hence, 
after passing to a subsequence, we have diffeomorphisms $F_k:\Mbar\rightarrow
\Mbar_k$ and a Riemannian metric $g$ on $\Mbar$ such that $F^*_kg_k\rightarrow
g$ in $C^r(\Mbar)$ for all $r<2$.  It remains to show that $g\in 
C^2_*(\Mbar)$.

Note that $F^*_k \Ric_{M_k}=\Ric_{F^*_kg_k}$ has uniformly bounded $L^\infty$
norm, when measured via $F^*_kg_k$, hence when measured via $g$.  Identities 
of the form (2.0.7)--(2.0.8) imply $F^*_k \Ric_{M_k}\rightarrow \Ric_M$
in $H^{-\ep,p}(M)$, for each $\ep>0,\ p<\infty$.  But the observation above
implies some subsequence converges $\text{weak}^*$ in $L^\infty$.
It follows that
$$
\Ric_M\in L^\infty(M).
$$
A similar argument gives $\Ric_{\pa M}\in L^\infty(\pa M)$, and also similarly
we obtain for the mean curvature $H$ of $\pa M\hookrightarrow \Mbar$ that
$H\in\Lip(\pa M)$.  Thus Theorem 2.1 applies to give $g\in C^2_*(\Mbar)$,
in boundary harmonic coordinates.  This finishes the proof of Theorem 3.1.

$\text{}$\newline
{\smc Remark}.  Invoking the definition of the Gromov-Hausdorff topology
(cf.~[Gr]) we can show that $\overline{\Cal{M}(R_0,i_0,S_0,d_0)}$ is compact
in the Gromov-Hausdorff topology and $C^r$-convergence is equivalent to
Gromov-Hausdorff convergence on this compact set, for any $r\in[1,2)$.

$\text{}$

\heading
\S{3.4}: Convergence of the geodesic flow and implications
\endheading

$\text{}$

Here we establish results on the limiting behavior of geodesic flows
under $C^1$ convergence, and implications for the injectivity radius,
that sharpen and generalize some of the results
of [Sak].  We work in the following setting.  Let $\Mbar$ be a fixed
compact manifold, with a $C^3_*$ coordinate system.  Let $g,g_k$ be metric 
tensors on $\Mbar$.  Assume
$$
g,\,g_k\in C^2_*(\Mbar),\quad g_k\rightarrow g\ \text{ in }\ C^1(\Mbar).
\tag{3.4.1}
$$
Say these metric tensors define exponential maps
$$
\Exp_p:T_p\Mbar\supset U\rightarrow \Mbar,\quad
\Exp_{k,p}:T_p\Mbar\supset U\rightarrow \Mbar,
\tag{3.4.2}
$$
where $U$ is a neighborhood of $0\in T_p\Mbar$ if $p\in M^{int}$ and $U$
consists of $v\in T_p\Mbar$ in a neighborhood of $0$ that point into 
$M^{int}$, transversally to $\pa M$, if $p\in\pa M$.  
Here is our first result.

\proclaim{Proposition 3.4.1} Under the hypotheses listed above,
given $v\in U\subset T_p\Mbar$,
$$
\lim\limits_{k\rightarrow\infty}\, \Exp_{k,p}(v)=\Exp_p(v).
\tag{3.4.3}
$$
\endproclaim
\demo{Proof} Let $X_k, X$ denote the vector fields on $T\Mbar$ generating 
these flows, essentially the Hamiltonian vector fields associated with 
$g_k$ and $g$.  Thus the coefficients of these vector fields belong to 
$C^1_*$ and there is $C^0$-convergence $X_k\rightarrow X$.  
Orbits of $X_k$ are
$$
y_k(t)=(\gamma_k(t),\gamma'_k(t)),
\tag{3.4.4}
$$
and we have $y'_k(t)=X_k(y_k(t))$ uniformly bounded, hence both
$\gamma_k$ and $\gamma'_k$ uniformly Lipschitz.  It follows from
Arzela's theorem that a subsequence converges locally uniformly:
$$
\gamma_{k_\nu}\rightarrow\sigma,\quad \gamma'_{k_\nu}\rightarrow \sigma'.
\tag{3.4.5}
$$
Clearly $(\sigma,\sigma')$ is an orbit of $X$, and $\sigma(0)=p,\ 
\sigma'(0)=v$.  Since (by Osgood's theorem) the flow generated by $X$
is unique, it follows that $\sigma\equiv\gamma$.  Since this is true for
any convergent subsequence, the result (3.4.3) follows.
\enddemo

We next discuss the injectivity radius.  We need to modify the definition
used for (3.0.2), since in the present case $\Exp_p:U\rightarrow\Mbar$
need not be $C^1$ or even Lipschitz.  If $p\in M^{int}$ and $\rho\le 
\dist(p,\pa M)$, we will say
$$
\wti(p,g)\ge\rho
\tag{3.4.6}
$$
provided that, for each unit vector $v\in T_pM^{int}$,
$$
|t|<\rho\Longrightarrow \dist(p,\Exp_p(tv))=|t|,
\tag{3.4.7}
$$
i.e., provided $\gamma_v(s)=\Exp_p(sv)$ is length-minimizing from $p$ to
$\gamma_v(t)$ for $|t|<\rho$.  This does imply injectivity:

\proclaim{Proposition 3.4.2} If $\wti(p,g)\ge\rho$, then
$$
\Exp_p:B_\rho(0)\longrightarrow M^{int}\ \text{ is one-to-one},
\tag{3.4.8}
$$
where $B_\rho(0)=\{v\in T_pM^{int}:g(v,v)<\rho^2\}$.
\endproclaim
\demo{Proof} If $v,w\in T_pM^{int}$ are unit vectors and 
$\gamma_v(s)=\gamma_w(t)$ with $|s|,|t|<\rho$, 
then the condition (3.4.7) forces $s=t$ (maybe after
changing the sign of $w$).  If $v\neq w$, these geodesics must intersect 
non-tangentially at $q=\gamma_v(t)=\gamma_w(t)$, i.e., at a positive angle.
Then a standard construction produces, for small $\ep>0$, a curve from $p$
to $\gamma_v(t+\ep)$ shorter than $|t|+\ep$, contradicting (3.4.7).
\enddemo

We next have the following semicontinuity result. 

\proclaim{Proposition 3.4.3} In the setting of Proposition 3.4.1,
$$
\wti(p,g_k)\ge \rho_0\ \forall\ k\Longrightarrow
\wti(p,g)\ge \rho_0.
\tag{3.4.9}
$$
\endproclaim
\demo{Proof} Given $v\in T_p M^{int}$ such that $g(v,v)=1$, we set
$v_k=v/\sqrt{g_k(v,v)}$, so $g_k(v_k,v_k)=1$ and $v_k\rightarrow v$.
We are given that for all $k$,
$$
|t|<\rho_0\Longrightarrow \dist_k(p,\Exp_{k,p}(tv_k))=|t|,
\tag{3.4.10}
$$
where $\dist_k$ denotes distance as determined by the metric tensor $g_k$.
Now it is elementary that
$$
\dist_k(p,q)\rightarrow \dist(p,q),\ \text{ as }k\rightarrow \infty,
\tag{3.4.11}
$$
uniformly for $q\in \Mbar$.  Hence
$$
\dist_k(p,\Exp_p(tv))\rightarrow \dist(p,\Exp_p(tv)),
\tag{3.4.12}
$$
while Proposition 3.4.1 together with (3.4.11)
implies
$$
\dist_k(p,\Exp_{k,p}(tv_k))-\dist_k(p,\Exp_p(tv))\rightarrow 0.
\tag{3.4.13}
$$
so we have
$$
|t|<\rho_0\Longrightarrow \dist(p,\Exp_p(tv))=|t|,
\tag{3.4.14}
$$
as desired.
\enddemo

We can also define the following sort of ``boundary injectivity radius.''
We say
$$
\wti_b(\Mbar,g)\ge\rho
\tag{3.4.15}
$$
provided that for each $p\in\pa M$, inward normal $\nu_p\in T_p\Mbar$
(orthogonal to $T_p(\pa M)$),
$$
0\le t<\rho\Longrightarrow \dist(\Exp_p(t\nu_p),\pa M)=t.
\tag{3.4.16}
$$
Parallel to Proposition 3.4.2, we have:

\proclaim{Proposition 3.4.4} If $\wti_b(\Mbar,g)\ge\rho$, then
$$
\Phi:\pa M\times[0,\rho)\longrightarrow \Mbar,
\tag{3.4.17}
$$
given by
$$
\Phi(p,t)=\Exp_p(t\nu_p),
\tag{3.4.18}
$$
is one-to-one.
\endproclaim

Then, parallel to Proposition 3.4.3, we have
\proclaim{Proposition 3.4.5} In the setting of Proposition 3.4.1,
$$
\wti_b(\Mbar,g_k)\ge\rho_0\ \forall\ k\Longrightarrow 
\wti_b(\Mbar,g)\ge\rho_0.
\tag{3.4.19}
$$
\endproclaim

The proofs of these results are very similar to those of their counterparts 
above.  Finally, we have the following significant implication for the 
geometric convergence obtained in Theorem 3.1.

\proclaim{Corollary 3.4.6} If $(\Mbar,g)$ is a limit of $(\Mbar_k,g_k)\in
\Cal{M}(R_0,i_0,S_0,d_0)$ as in Theorem 3.1, then
$$
\aligned
\wti(p,g)&\ge \min(i_0,\dist(p,\pa M)),
\quad \forall\ p\in M^{int}, \\
\wti(p,g|_{\pa M})&\ge i_0,\quad \forall\ p\in\pa M, \\
\wti_b(\Mbar,g)&\ge i_0.
\endaligned
\tag{3.4.20}
$$
\endproclaim

$$\text{}$$
{\bf 4. Gel'fand inverse boundary problem}
\newline {}\newline

In this section we prove uniqueness and stability
and provide a reconstruction procedure for the 
inverse boundary spectral problem. To fix notations, assume that 
$(\Mbar,g,\p M)$ is a compact, connected manifold, with nonempty 
boundary, provided with a metric tensor $g$ with some limited
smoothness (specified more precisely below).
Let $\Delta^N$ be the Neumann Laplacian.  Denote by 
$$
0=\lambda_1 < \lambda _2 \leq \cdots
$$
the eigenvalues 
(counting multiplicity) of $-\Delta^N$ and
by  
$$
\phi_1 ={\Vol}(M)^{-1/2},\ \phi_2, \dots,
$$
the corresponding orthonormalized eigenfunctions.

The {\it Gel'fand inverse boundary problem} is the problem of the
reconstruction of $(\Mbar,g)$ from its boundary spectral data,
i.e., the collection
$(\p M,\{\lambda _k, \phi _k|_{\p M}\}_{k=1}^{\infty})$.

$\text{}$ \newline
{\smc Remark}.  The data used in the original formulation of the
Gel'fand inverse boundary problem [Ge] 
consists of the trace on $\p M$ of the resolvent
kernel, $R_{\lambda}(x,y)$ of $\Delta^N$ given for all 
$\lambda \in {\CC}\setminus\text{spec}\,\Delta^N,\ x,y \in \p M$. 
The equivalence of these data and boundary spectral 
data, and also dynamic inverse boundary
data for the corresponding wave, heat and non-stationary Schr{\"o}dinger
equations is proven in [KKL], [KLM], as are analogous results for the 
Dirichlet problem.
\newline $\text{}$

The main results of this section are the uniqueness result below and 
stability result in \S4.3.

\proclaim{Theorem 4.1}   Let $\Mbar$ be a compact, connected manifold with 
nonempty boundary and $C^2_*$ metric tensor.  Then the boundary
spectral data $(\p M, \, \{\lambda _k, \phi _k|_{\p M}\}_{k=1}^{\infty})$
determine the manifold $\Mbar$ and its metric $g$ uniquely.
\endproclaim

Such a result was established in the $C^{\infty}$ case in [BK1],
taking into account [Ta]; see also [KKL].
Our proof here will incorporate techniques from these papers, 
plus some additional arguments necessary to handle the reduced 
smoothness.  One essential role played by the $C^2_*$-hypothesis 
is that this implies non-branching of geodesics on $\Mbar$, including
geodesics passing transversally from $\pa M$.  In \S{4.1} we show that
$\Mbar$ is uniquely determined as a topological space.  We proceed in 
\S{4.2} to show that the differential structure and metric tensor on 
$\Mbar$ are uniquely determined.  In \S{4.3} we apply Theorem 4.1 
together with Theorem 3.1 to establish a result on the conditional stability
of this inverse boundary problem, building on results of [K2L].

$\text{}$

\heading
\S{4.1}: Determining the domain
\endheading

$\text{}$

We start with the introduction of some useful geometric objects.
Let $\Gamma \subset \p M $ be open and take $t\geq 0$. Then we set
$$
M(\Gamma,t) = \{x \in M: d(x,\Gamma) \leq t \},
\tag{4.1.1}
$$
the domain of influence of  $\Gamma$ at ``time'' $t$. This
terminology refers to the corresponding wave equation where
$M(\Gamma,t)$ is the subdomain of $M$ filled  by the time $t$ with
waves sent from $\Gamma$.

Now let $\underline{\Gamma}$ consist of a finite number
$\Gamma_1, \dots,\Gamma_m$ of subsets  $\Gamma$ and
$\underline{t}^+, \ \underline{t}^-$ be two $m$-dimensional vectors
with positive entries, 
$\underline{t}^+ =(t^+_1, \dots,t^+_m)$,
$\underline{t}^- =(t^-_1, \dots,t^-_m)$.
Then set
$$
M(\underline{  \Gamma}, \,\underline{t}^+, \, \underline{t}^-)=
\bigcap_{i=1}^m (M(\Gamma_i ,t^+_i) \setminus M(\Gamma_i ,t^-_i))
\subset M,
\tag{4.1.2}
$$
and define 
$$
{\bL}(\underline{  \Gamma}, \,\underline{t}^+, \, \underline{t}^-)=
{\F}L^2(M(\underline{  \Gamma}, \,\underline{t}^+, \, \underline{t}^-))
\subset \ell^2.
\tag{4.1.3}
$$
Here, ${\F}$ stands for the Fourier transform of functions from
$L^2(M)$,
$$
{\F}(u) = \{u_k\}_{k=1}^ {\infty}\in \ell^2,
\quad
u(x) = \sum_{k=1}^\infty u_k \phi_k(x), \quad  u_k=(u,\phi_k)_{L^2(M)},
\tag{4.1.4}
$$
and the subspace
$L^2(M(\underline{\Gamma}, \,\underline{t}^+, \, \underline{t}^-))$ 
consists of all functions in $L^2(M)$ with support in the set
$M(\underline{\Gamma}, \,\underline{t}^+, \, \underline{t}^-)$.

Two basic ingredients for the reconstruction of the manifold $\Mbar$
are the approximate controllability and Blagovestchenskii's formula.

The controllability result is an implication of Tataru's
unique continuation result for the wave
equation ([Ta], see also [Ho], [Ta2]). To describe it,
consider the wave equation
$$
\aligned
&(\p_t^2-\Delta  )u^f(x,t)=0\quad\hbox{in }M\times \RR_+ \\
&u^f|_{t=0}=0,\quad u_t^f|_{t=0}=0,\quad Nu^f|_{\p M\times \RR_+ }=f
\in C^1_0(\Gamma\times(0,T)),
\endaligned
\tag{4.1.5}
$$
where $N$ is the exterior unit normal field to $\p M$.
Using Tataru's theorem, it was shown in Theorem 3.10 of [KKL] that 
the following holds.

\proclaim{Proposition 4.1.1} For each $T>0$, the set
$\{u^f(T):f\in L^2(\Gamma\times(0,T))\}$ is a
dense subspace of $L^2(M(\Gamma,T))$.  
\endproclaim

(Actually Theorem 3.10 of [KKL] is written to address the Dirichlet boundary
condition, but the same argument works for the Neumann boundary condition.)

Blagovestchenskii's formula gives the Fourier coefficients $u_k^{f}(t)$
of a wave $u^{f}(\cdotp,t)$ in terms of the boundary spectral data,
$$
u_k^{f}(t) =
\int_0^t \int\limits_{\p M} f(x,t')\,
\frac {\sin {\sqrt {\lambda _k} (t-t' )}}{\sqrt {\lambda _k}}
\phi_k(x)\,dS_g\,dt'.
\tag{4.1.6}
$$
To prove this one starts with $\pa_t^2(u(t),\phi_k)_{L^2}=(\Delta u(t),
\phi_k)_{L^2}$ and applies Green's formula to get an inhomogeneous ODE
for $(u(t),\phi_k)_{L^2}$, yielding (4.1.6).

Note that in the formula (4.1.6) there appears the  Riemannian volume 
$dS_g$ of $\p M$, which we are not given.  However, we are given 
$\pa M$ as a $C^2$ manifold, so an arbitrarily chosen volume element
element has the form $dS=\kappa\, dS_g$, where $\kappa$ is $C^1$-smooth and
strictly positive.  We can construct the Fourier coefficients of the wave
$u^{\kappa f}(x,t)$ for any boundary source $f$, and despite our lack
of knowledge of $\kappa$, we do have the following.

\proclaim{Corollary 4.1.2}  Given $\Gamma\subset \pa M$ and $t>0$,
the boundary spectral data determine the subspace
$$
{\bL}(\Gamma,t) = \Cal{F}L^2(M( \Gamma, \,t))\subset\ell^2.
\tag{4.1.7}
$$
In fact, let $\{f_\nu:\nu\in\ZZ^+\}$ have dense linear span in $L^2(\Gamma
\times(0,T))$.  Then ${\bL}(\Gamma,t)$ is the closed linear span in $\ell^2$
of $\{\varphi_\nu:\nu\in\ZZ\}$, where $\varphi_\nu\in\ell^2$ is given by
$\varphi_{\nu,k}=u^{f_\nu}_k(t)$.

Thus we can find the orthoprojection
$P: \ell^2 \to {\bL}(\Gamma,t) $ to this subspace. 
\endproclaim

From here, using the elementary identities
$$
L^2\Bigl(\bigcap\limits_i S_i\Bigr)=\bigcap\limits_i L^2(S_i),\quad
L^2(A_i\setminus B_i)=L^2(A_i)\cap L^2(B_i)^\perp,
\tag{4.1.8}
$$
we deduce that the boundary spectral data uniquely determine the subspaces
${\bL}(\underline{  \Gamma}, \,\underline{t}^+,$  
$\underline{t}^-)$ of $\ell^2$,
for any $\underline{  \Gamma}, \,\underline{t}^+, \, \underline{t}^-$
with arbitrary $m$ (compare [KKL] and [Be1]).
In particular, for any $\underline{  \Gamma}, \,\underline{t}^+,
\, \underline{t}^-$ we can see if
${\bL}(\underline{  \Gamma}, \,\underline{t}^+, \, \underline{t}^-) 
= \{0\}$ or not. Equivalently, we can see if
$M(\underline{  \Gamma}, \,\underline{t}^+, \, \underline{t}^-)$
contains an open ball or not.

Next, let $h\in C(\p M)$. We can ask if $h$ is the boundary distance
function for some $x \in M$.  To this end, we choose points
$z_j\in \p M$, $j=1,\dots,m$, their small neighborhoods
$\Gamma_1,\dots,\Gamma_m$ and numbers $t^\pm_j=h(z_j)\pm 1/m$. 
When $m\to \infty$, the fact that ${\bL}(\underline{\Gamma},
\,\underline{t}^+, \, \underline{t}^-) \neq \{0\} $ for any $m$
determines whether
there is a point $x\in M$ such that $h(z)=\text{dist}(z,x),  \, z \in \p M $.
Thus we have shown that the boundary spectral data
determine the image in $L^{\infty}(\p M)$ of the 
boundary distance representation 
$R$. Here, $R:\Mbar \to C(\p M)$ is defined by 
$$
R(x) = r_x(\cdot), \quad r_x(z) = \text{dist}(x,z), \quad z \in \p M.
\tag{4.1.9}
$$
(Compare [KKL] and [Ku1]).
Clearly, the map $R$ is Lipschitz continuous. Moreover, under the assumptions
of Theorem 4.1 it is injective. To see this, let $r_x = r_y$ and let
$z \in \p M$ be a point of minimum of these functions. Then both $x$ and $y$
lie on the normal geodesic to $ \p M$ starting in $z$ at the same arclength
$r_x(z)= r_y(z)$. 
As the metric $g \in C^2_*(\Mbar)$, it follows from Corollary 2.5.2
that this normal geodesic does not branch. Therefore, $x=y$.

Since $\Mbar$ is compact, injectivity and continuity imply that $R$ is a
homeomorphism, i.e., $R(\Mbar)$ with the distance inherited from 
$L^{\infty}(\p M)$ and $(\Mbar,g)$ are homeomorphic, and thus
$R(\Mbar)$ can be identified with $\Mbar$ as a topological manifold. 
We have established the following.

\proclaim{Proposition 4.1.3} Assume $(\Mbar_1,g_1)$ and $(\Mbar_2,g_2)$ 
satisfy the hypotheses of Theorem 4.1.  If they have identical boundary 
spectral data, including $\pa M_1=\pa M_2=X$, as $C^2$ manifolds,
then there is a natural correspondence of $R(\Mbar_1)$ and 
$R(\Mbar_2)\subset C(X)$, producing a uniquely defined homeomorphism
$$
\chi:\Mbar_1\longrightarrow \Mbar_2.
\tag{4.1.10}
$$
\endproclaim

$\text{}$

\heading
\S{4.2}: Determining the metric
\endheading

$\text{}$

Our next goal is to reconstruct the differential and Riemannian
structures on $M$. To this end, let us return to 
${\bL}( \underline{  \Gamma}, \,\underline{t}^+, \, \underline{t}^-)$ 
and consider the orthoprojection 
${P}( \underline{  \Gamma}, \,\underline{t}^+, \, \underline{t}^-)$
of $\ell_2$ onto this subspace. Then,
$$
\left({P}( \underline{  \Gamma}, \,\underline{t}^+, \, \underline{t}^-)e_i, 
e_j\right)_{\ell^2}=
\int_{M( \underline{  \Gamma}, \,\underline{t}^+, \, \underline{t}^-)}
\phi _i(x) \phi _j(x) \,dV_x,
\tag{4.2.1}
$$
where $e_j=(0,\dots,0,1,0,\dots)$ is the sequence having a 1 at the $j$th 
place.  Also,
$$
\int_{M( \underline{  \Gamma}, \,\underline{t}^+, \, \underline{t}^-)}
\phi _1(x)^2 \,dV_x =
\frac{\Vol(M( \underline{  \Gamma}, \,\underline{t}^+, \, 
\underline{t}^-))}{\Vol(M)}.
\tag{4.2.2}
$$

Next choose a sequence $(\underline{\Gamma}_k, \,\underline{t}^+_k, \, 
\underline{t}^-_k),\ k=1,2,\dots$,
with $m_k \rightarrow \infty$, where $m_k$ is
the dimension of $\underline{t}^{\pm}_k,$ so that 
$M( \underline{  \Gamma}_k, \,\underline{t}^+_k, \, \underline{t}^-_k)$
shrinks to $\{x\}$ when $k\to \infty$.
Then by formulae (4.2.1)--(4.2.2) we see that
$$
\lim_{k\to \infty}\left({P}( \underline{  \Gamma}_k, \,\underline{t}^+_k, 
\, \underline{t}^-_k)e_1, \,e_j\right)
\cdot \left({P}( \underline{  \Gamma}_k, \,\underline{t}^+_k, 
\, \underline{t}^-_k)e_1, \,e_1\right)^{-1/2}=\phi_j(x).
\tag{4.2.3}
$$
Thus we can find values of the eigenfunctions 
$\phi _k(x)$ for all $k=1,2, \dots$ and $x \in M$.

To proceed further we need an auxiliary statement about the
properties of the eigenfunctions. Let
$\Phi$ be  the space of all
finite linear combinations of $\phi _k, \ k=1,2,\dots $.

\proclaim{Lemma 4.2.1}
Under the assumptions of Theorem 4.1, $\Phi$ is dense in 
the space $\{u\in H^{s,p}(M): N u|_{\pa M}=0\}$
for any  $s\in[0,3),\ p\in(1,\infty)$.
Moreover, if $x \in M^{{int}}$, there are
$n$ indices $k(1),\dots,k(n)$ (depending on $x$) 
and a neighborhood $U$ of $x$
such that $\phi _{k(1)}(x), \dots,\phi _{k(n)}(x)$ form a $C^3_*$-smooth
coordinate system in $U$.
\endproclaim
\demo{Proof}
Assuming that $g \in C^{r}, \, r \in (1,2)$, consider,
for any $p\in (1,\infty)$,   the Neumann 
Laplacian,  $\Delta^N_p$, with domain
$$
{\Cal  D}(\Delta^N_p)=\{u\in H^{2,p}(M): N u|_{\pa M}=0\}.
\tag{4.2.4}
$$
Denote by $ e^{t\Delta^N_p}, \, t \geq 0$  the
corresponding contraction semigroup and by   $(-\Delta^N_p)^s,\ s>0$, the
real powers  of $-\Delta^N_p$, defined for $s\in(0,1)$
via subordination. By  Stein's 
Littlewood-Paley theory for symmetric diffusion semigroups (cf.~[St]),
$$
Y^{s,p}={\Cal  D}((-\Delta^N_p)^{s/2}),\quad s\in [0,\infty),
\tag{4.2.5}
$$
is a complex interpolation scale, in $s$, for each $p\in (1,\infty)$.
In particular,
$$
{\Cal  D}((-\Delta^N_p)^{s/2})=H^{s,p}(M), \quad 
\hbox{for} \, \,0 \le s \le 1.
\tag{4.2.6}
$$
Hence, for $0\le s\le 1$,
$$
\aligned
{\Cal  D}((-\Delta^N_p)^{1+s/2})& =  \{ u\in{\Cal D}(-\Delta^N_p): 
\Delta u\in{\Cal D}((\Delta^N_p)^{s/2})\} \\ 
&=\{u\in H^{2,p}(M): N u|_{\pa M}=0, 
\hbox{ and }\Delta u\in H^{s,p}(M)\}.
\endaligned
\tag{4.2.7}
$$
If we require
$g\in C^2_*(M),$
we can use regularity results to obtain that
$$
{\Cal D}((-\Delta^N_p)^{1+s/2})=
\{u\in H^{2+s,p}(M): N u|_{\pa M}=0\},\quad 0\le 
s<1.
\tag{4.2.8}
$$
Since $(1-\Delta^N_p)^{s/2}:Y^{s,p}\rightarrow L^p(M)$ is an
isomorphism, acting bijectively on $\Phi$,  
the desired density of $\Phi$ will follow from the 
density of $\Phi$  in $L^p(M)$.  We now demonstrate this density.

Suppose $f\in L^{q}(M)$ (with $q=p'$)
and $\langle f,u\rangle=0$ for all $u\in \Phi$.
If $q \ge 2$ then $f\in L^2(M)$ and clearly $f\equiv 0$.  So we need only 
worry about the case $q<2$.  
Note that $\langle f,u\rangle=0$ implies 
$$
\langle e^{t\Delta^N}f,u\rangle =0,\quad \forall\ t\ge 0,\ u\in \Phi.
\tag{4.2.9}
$$
Now $e^{t\Delta^N_q}$ is a holomorphic semigroup on $L^{q}(M)$, so 
for all $t>0$,
$$
e^{t\Delta^N_q }f\in{\Cal D}(\Delta^N_q)\subset H^{2,q}(M)\subset 
L^{q_2}(M),
\tag{4.2.10}
$$
with $q_2>q$, by the Sobolev embedding theorem.  Iterating this and 
using the semigroup property gives
$$
f\in L^{q}(M) \Longrightarrow e^{t\Delta^N}f\in L^2(M) \quad 
\forall\  t>0.
\tag{4.2.11}
$$
Hence (4.2.9) implies $e^{t\Delta^N}f=0$, for all $t>0$.  But 
$e^{t\Delta^N}f\rightarrow f$
in $L^{q}(M)$ as $t\searrow 0$, so $f=0$.  This completes the proof
of the first statement of the lemma.

To demonstrate the second statement,
we first note that since $s<3$ and $p< \infty$
are arbitrary, it follows from the first part of the lemma that
$$
C^{2}_0(M^{int}) \subset
\text{closure of }\Phi \text{ in } C^2(\Mbar).
\tag{4.2.12}
$$
Let now $x \in M^{{int}}$ and 
$(x^1, \cdots, x^n)$ be some local coordinates near $x$. 
Denote by  $T_x: \Phi \to \RR^n$ the map,
$$
T_x(u) = (\p_1 u(x), \cdots, \p_nu(x)).
\tag{4.2.13}
$$
It follows from (4.2.12) that $T_x(\Phi) = \RR^n$, i.e., there are
indices $k(1),\dots,k(n)$ (depending on $x$) such that 
$\nabla \phi_{k(i)}(x),\ i=1, \dots, n$, are linearly independent. 
Moreover, we know that the eigenfunctions $\phi_k \in C^3_*( M^{{int}})$.
This proves the second statement of the lemma.
\enddemo

Having this, we are in a position to refine our statement about the 
homeomorphism $\chi:\Mbar_1\rightarrow\Mbar_2$ established in Proposition 
4.1.3.
\proclaim{Proposition 4.2.2} Let $(\Mbar_1,g_1)$ and $(\Mbar_2,g_2)$ be as
in Proposition 4.1.3, with identical boundary spectral data.  Then the 
map $\chi$ in (4.1.10) has the property that
$$
\chi:M_1^{int}\longrightarrow M_2^{int}\ \text{ is a }
C^2\text{-diffeomorphism, and }\ \chi^*g_2=g_1.
\tag{4.2.14}
$$
\endproclaim
\demo{Proof} We use the fact that, if $\{\phi_j\}$ are the normalized
eigenfunctions of $\Delta^N$ on $M_1$ and $\{\tilde{\phi}_j\}$ those on 
$M_2$, then, as a consequence of (4.2.3),
$$
\phi_j(x)=\tilde{\phi}_j(\chi(x)).
\tag{4.2.15}
$$
Given $p\in M_1^{int}$, there exist indices $k(1),\dots,k(n)$ such that
$\tilde{\phi}_{k(1)},\dots,\tilde{\phi}_{k(n)}$ form a local coordinate
system on a neighborhood of $\tilde{p}=\chi(p)$.  Then, for $x$ near $p$,
we have
$$
x\mapsto (\phi_{k(1)}(x),\dots,\phi_{k(n)}(x))=
\bigl(\tilde{\phi}_{k(1)}(\chi(x)),\dots,\tilde{\phi}_{k(n)}(\chi(x))\bigr) 
\mapsto \chi(x)
\tag{4.2.16}
$$
a composition of $C^2$ smooth maps.  Thus $\chi$ in (4.2.14) is $C^2$ smooth 
on a neighborhood of each $p\in M_1^{int}$, hence on $M_1^{int}$.  
Interchanging the roles of $\Mbar_1$ and $\Mbar_2$, we have the same 
result for $\chi^{-1}$.

Finally, we show that the metric tensor is uniquely determined.
For notational simplicity, just consider $M=M_1$.
Let $(x^1,\dots,x^n)$ be a $ C^3_* \,$
coordinate system in a domain 
$U \subset M^{{int}}$, e.g., the one obtained earlier from 
the eigenfunctions. Then, for all $k=1,2, \dots,$
$$
-g^{ij}(x) \p_i \p_j \phi _k(x) -b^i(x) \p _i \phi_k(x) = \lambda _k
\phi _k(x), \quad b^i = g^{-1/2} \p _j(g^{1/2}g^{ij}), 
\tag{4.2.17}
$$
where all
eigenfunctions $\phi_k$ and, henceforth, their derivatives as well as
$\lambda _k$ are already found.  Let us consider equations
(4.2.17) as linear equations for $g^{ij}(x)=g^{ji}(x), \ b^i(x).$
Using again (4.2.12), we see that the map $\widetilde{ T }_x:
\Phi \to
\RR^{ n+n(n+1)/2} $, 
$$
\widetilde{ T} _x(u) =(\p_i u(x), \, 
\p_i \p_j u(x): \, i \leq j =1, \cdots, m),
\tag{4.2.18}
$$
is surjective (compare with (4.2.13)).  Thus equations
(4.2.17) are uniquely solvable since the 
2-jets of the eigenfunctions $\phi _k$ at $x$ span the whole space
$\RR^{n+n(n+1)/2}$.

It follows that the diffeomorphism $\chi$ in (4.2.14) pulls the metric tensor
$g_2$ back to $g_1$.  The proof of Proposition 4.2.2, and hence of 
Theorem 4.1, is complete.
\enddemo

We can get some more insight into how the geometry of $(\Mbar,g)$ is
determined by the boundary spectral data, particularly through (4.2.3),
by examining further the maps
$$
\Psi_k:\Mbar\longrightarrow \RR^k,\quad \Psi_k(x)=(\phi_1(x),\dots,\phi_k(x)).
\tag{4.2.19}
$$
For simplicity we assume the eigenfunctions are arranged to be real valued.
The argument proving Proposition 4.2.3 shows that for each compact $K\subset
M^{int}$, there exists $k$ such that $\Psi_k$ restricted to $K$ is an 
embedding.  In fact, we can do better than that, though for no $k$ will
$\Psi_k$ be an embedding of $\Mbar$, since $D\Psi_k(x)$ annihilates the
normal to $\pa M$ for each $x\in\pa M$, each $k$.  Note that Lemma 4.2.1
implies the space of restrictions of elements of $\Phi$ to $\pa M$ is dense
in $C^2(\pa M)$, so we can find $k_0=k_0(\Mbar,g)$ such that
$$
\Psi_k:\pa M\rightarrow \RR^k\ \text{ is an embedding, for }\ k\ge k_0.
\tag{4.2.20}
$$
We now augment $\Psi_k$ to
$$
\Psi^\#_k:\Mbar\longrightarrow \RR^{k+1},\quad
\Psi^\#_k(x)=(\psi_0(x),\phi_1(x),\dots,\phi_k(x)),
\tag{4.2.21}
$$
where $\psi_0\in C^3_*(\Mbar)$ is the eigenfunction for the {\it Dirichlet
problem}, with smallest eigenvalue:
$$
\Delta\psi_0=-\mu_0 \psi_0,\quad \psi_0\bigr|_{\pa M}=0,
\tag{4.2.22}
$$
normalized by
$$
\int\limits_M |\psi_0(x)|^2\, dV=1,
\tag{4.2.23}
$$
and let us insist $\psi_0(x)>0$ on $M^{int}$.  Hopf's principle implies
$N\psi_0(x)\neq 0,\ \forall\ x\in\pa M$.  Hence the map (4.2.21) has the
property
$$
D\Psi^\#_k(x)\ \text{ is injective},\quad \forall\ x\in\pa M,\ k\ge k_0.
\tag{4.2.24}
$$
Thus injectivity holds on a collar neighborhood of $\pa M$ in $\Mbar$.
Combined with our previous observation about $\Psi_k$ embedding compact
$K\subset M^{int}$, this implies that there exists $k_1=k_1(\Mbar,g)$
such that
$$
D\Psi^\#_k(x)\ \text{ is injective},\quad \forall\ x\in\Mbar,\ k\ge k_1.
\tag{4.2.25}
$$
Also we have $\Psi^\#_k$ embedding $\pa M$.  Perhaps increasing $k_1$ to be 
sure points are completely separated, we have the following result.

\proclaim{Proposition 4.2.3} Let $\Mbar$ be as in Theorem 4.2.  Then there
exists $k_2=k_2(\Mbar,g)$ such that
$$
\Psi^\#_k:\Mbar\longrightarrow \RR^{k+1}\ \text{ is a }
C^2\text{-embedding}, \quad \forall\ k\ge k_2.
\tag{4.2.26}
$$
\endproclaim

Having this, we can improve Proposition 4.2.2, as follows.
\proclaim{Corollary 4.2.4} In the setting of Proposition 4.2.2, we have
$$
\chi:\Mbar_1\longrightarrow \Mbar_2\ \text{ is a }C^2\text{-diffeomorphism}.
\tag{4.2.27}
$$
\endproclaim
\demo{Proof} Apply the construction (4.2.19)--(4.2.26) to obtain embeddings
$$
\Psi^\#_k:\Mbar_1\longrightarrow \RR^{k+1},\quad
\widetilde{\Psi}^\#_k:\Mbar_2\longrightarrow\RR^{k+1},
\tag{4.2.28}
$$
for $k$ sufficiently large.  Here we take $\Psi^\#_k$ as in (4.2.21),
with $(\Mbar,g)=(\Mbar_1,g_1)$, and define
$$
\widetilde{\Psi}^\#_k(x)=(\tilde{\psi}_0(x),\tilde{\phi}_1(x),\dots,
\tilde{\phi}_k(x))
$$
on $(\Mbar_2,g_2)$ in the analogous fashion.  Using Proposition 4.2.2, we
see that $\chi$ pulls back the volume element of $(M_2,g_2)$ to that of
$(M_1,g_1)$, and deduce that
$$
\psi_0(x)=\tilde{\psi}_0(\chi(x)).
\tag{4.2.29}
$$
In concert with (4.2.15), this gives
$$
\Psi^\#_k=\widetilde{\Psi}^\#_k\circ \chi.
\tag{4.2.30}
$$
Since both maps in (4.2.28) are $C^2$-embeddings (onto the same range),
the result (4.2.27) follows from the implicit function theorem.
\enddemo

$\text{}$ \newline
{\smc Remark 1}.
One can readily strengthen the regularity results on $\chi$, given above,
to regularity of class $C^3_*$.
Details can be left to the reader.

$\text{}$ \newline
{\smc Remark 2}.
The reason we do not use in this section the construction of [KKL]
to recover the 
differential and Riemannian structure of $(M,g)$ is the following. 
In [KKL] we use as coordinates some distance functions on $M$.
However, it is well known (e.g., [DTK]) that the 
resulting coordinates, 
in principle, lose two orders of regularity, i.e., in our case are
just $C^1_*$ regular. Clearly, it is rather inconvenient to work with such
coordinates.

$\text{}$ \newline
{\smc Remark 3}.
There is an analogue of Theorem 4.1
in the case of the Dirichlet boundary spectral data
$(\p M, \{\lambda _k,$ $ N\psi _1|_{\p M}\}_{k=1}^{\infty})$.
The proof is actually simpler. However, it requires some modifications
because we can no longer find 
$ \Vol(M( \underline{  \Gamma}, \,\underline{t}^+, \, \underline{t}^-)) $.
Therefore, instead of $\psi _k(x)$ we find
$$
\xi_k(x) = \frac{\psi _k(x)}{\psi _1(x)}, \quad x \in M^{{int}}, \quad
k >1.
$$
However, an analog of Lemma 4.2.1 is valid and we can construct
$C^3_*$-smooth coordinates in the vicinity of any $x \in M^{{int}}$
using functions $\xi_k, \, k>1$. Moreover, the map
$\tilde{ T }_x $, given by (4.2.18) is surjective on the linear span
of the functions $\xi_k$. Since these functions satisfy equations
$$
-g^{ij}(x) \p_i \p_j \xi  _k(x) -
\tilde{b^i}(x)
\p _i \xi _k(x) =(\lambda _k-\lambda _1)
\xi  _k(x), 
$$
with the same metric tensor $g^{ij}$ but different $\tilde{b^i}$, we 
can use these equations to find $g^{ij}(x),  \, x \in M^{{int}}$.
Having found the metric inside $M$ we can return to the 
analogs of equations (4.2.1), (4.2.2)
and reconstruct $\psi _k(x), \, x \in M^{{int}}$.

$\text{}$

\heading
\S{4.3}: Stabilization of the inverse problem
\endheading

$\text{}$

In this section we consider stabilization of inverse problems using 
geometric convergence results and apply them to the Gel'fand problem.
The basic thrust of our argument provides an illustration of a general
``stabilization principle for inverse problems,'' which we can describe
abstractly as follows.

Suppose $\Cal{M}$ is a collection of objects, one element $M$ of which
you want to identify via the observation of data $\Cal{D}(M)$, in some
set $\Cal{B}$ of observable data.  Suppose $\Cal{M}$ and $\Cal{B}$
have natural topologies, and the map 
$$
\Cal{D}:\Cal{M}\longrightarrow \Cal{B}
$$
has been shown to be continuous.  Suppose the uniqueness problem has
been solved, so you know this map is one-to-one.  However, typically
such a map does not have a continuous inverse.  This is a standard 
situation in the study of inverse problems, giving rise to the
phenomenon of ill-posedness.  This problem is made more acute by the
fact that what one measures is not exactly equal to $\Cal{D}(M)$,
but only an approximation to it.

The key to the stabilization, which is useful for a wide variety of 
inverse problems, requires an {\it a priori} knowledge that the object
$M$ one wants to identify actually belongs to a subset $\Cal{M}_0$ of
$\Cal{M}$, and that furthermore its closure in $\Cal{M},\ 
\overline{\Cal{M}}_0$, is {\it compact}.  In that case the restriction
$$
\Cal{D}:\overline{\Cal{M}}_0\longrightarrow \Cal{B}
$$
is automatically a {\it homeomorphism} of $\overline{\Cal{M}}_0$
onto its range in $\Cal{B}$.  Thus, when trying to identify the desired
object $M$, one minimizes some measure of the difference between the 
calculated data $\Cal{D}(M_j)$ and the observed data, while constraining
$M_j$ to belong to $\overline{\Cal{M}}_0$.

Having set up the abstract stabilization principle, we show how it applies
to the Gel'fand problem.

To prepare for this discussion, let us set up some notation.
Denote by $\Cal{M}_X(C^2_*)$ the set of compact, connected manifolds
$\Mbar$ with nonempty boundary $X$, endowed with a metric tensor in
$C^2_*(\Mbar)$.  Given $(\Mbar,g)\in\Cal{M}_X(C^2_*)$, set
$$
\Cal{D}(\Mbar,g)=\{\lambda_j,\phi_j|_{X}\}_{j=1}^\infty,
\tag{4.3.1}
$$
the right side denoting the boundary spectral data of $(\Mbar,g)$.
We have 
$$
\Cal{D}:\Cal{M}_X(C^2_*)\longrightarrow \Cal{B}_X,
\tag{4.3.2}
$$
where $\Cal{B}_X$ denotes the set of sequences $\{\mu_j,\psi_j:j\ge 1\}$,
with $\mu_j\in\RR^+,\ \mu_j\nearrow+\infty$, and $\psi_j\in L^2(X)$,
modulo an equivalence relation, which can be described as follows.
We say $\{\mu_j,\psi_j\}\sim \{\mu_j,\tilde{\psi}_j\}$ 
if $\psi_j(x)=\alpha_j \tilde{\psi}_j(x)$ 
for some $\alpha_j\in\CC, |\alpha_j|=1$.
More generally, if $\mu_{k_0}=\cdots=\mu_{k_1}$, we allow
$$
\psi_j(x)=\sum\limits_{k=k_0}^{k_1} \alpha_{jk} \tilde{\psi}_k(x),
\quad j = k_0, \cdots, k_1,
\tag{4.3.3}
$$
for a unitary $l\times l$ matrix $(\alpha_{jk}),\ l=k_1-k_0+1$.
The content of Theorem 4.1 is that the map (4.3.2) is one-to-one.

There are natural topologies one can put on the sets in (4.3.2).  On 
$\Cal{M}_X(C^2_*)$ one has the topology of $C^{r}$ 
convergence,
for any $r\in (1,2)$, defined in \S{3}.
On $\Cal{B}_X$ one has a topology described as follows.
We describe when $\{\mu_j^\nu,\psi_j^\nu\}\rightarrow \{\mu_j,\psi_j\}$,
as $\nu\rightarrow\infty$.  First we require $\mu_j^\nu\rightarrow \mu_j$
for each $j$.  Next, if $\mu_k$ is simple, i.e., different from $\mu_{k-1}$
and $\mu_{k+1}$, we require $\alpha^\nu\psi^\nu_k\rightarrow \psi_j$
in $L^2(X)$, for some $\alpha^\nu\in\CC,\ |\alpha^\nu|=1$.
More generally, if $\mu_k$ has multiplicity $\ell$, say $\mu_k=\cdots=
\mu_{k+\ell-1}$, we require that there exist unitary $\ell\times\ell$
matrices $(\alpha^\nu_{ij})_{k\le i,j\le k+\ell-1}$ such that
$$
\sum\limits_{j=k}^{k+\ell-1} \alpha^\nu_{ij} \psi^\nu_j
\longrightarrow \psi_i\ \text{ in }\ L^2(X).
\tag{4.3.4}
$$
Compare [K2L].

Given these topologies, it follows from standard techniques of perturbation
theory (cf.~[K]) that $\Cal{D}$ is continuous in (4.3.2).

Now the map (4.3.2) is  by no means invertible.  This is a standard situation 
encountered in the study of inverse problems, giving rise to the 
phenomenon of ill-posedness.  One wants to ``stabilize'' the inverse
problem, showing that certain {\it a priori} hypotheses on the domain
$(\Mbar,g)$ put it in a subset $K\subset \Cal{M}_X(C^2_*)$ having the property
that $\Cal{D}^{-1}$ can be shown to act continuously on the image of $K$.
The results of \S{3} provide a tool to accomplish this.

Recall the class $\Cal{M}(R_0,i_0,S_0,d_0)$ defined in \S{3}.  Given a 
boundary $X$, let $\Cal{M}_X(R_0,$ $i_0,S_0,d_0)$ denote the set of such 
manifolds with boundary $X$.  It follows from Theorem 3.1 that
$\overline{\Cal{M}_X(R_0,i_0,S_0,d_0)}$ is compact in the $C^{r}$
topology, for any $r\in (1,2)$, and is contained in $\Cal{M}_X(C^2_*)$.
We hence give $\overline{\Cal{M}_X(R_0,i_0,S_0,d_0)}$ the $C^{r}$
topology, and we see this is independent of $r$, for $r\in (1,2)$.

Combined with Theorem 4.1, these observations yield the following
conditional stability of the Gel'fand inverse problem.

\proclaim{Theorem 4.3.1} Given $R_0,i_0,S_0,d_0\in (0,\infty)$,
$$
\Cal{D}:\overline{\Cal{M}_X(R_0,i_0,S_0,d_0)}\longrightarrow \Cal{B}_X
$$
is a homeomorphism of $\overline{\Cal{M}_X(R_0,i_0,S_0,d_0)}$
onto its range, $\Cal{B}_X(R_0,i_0,S_0,d_0)$; hence
$$
\Cal{D}^{-1}:\Cal{B}_X(R_0,i_0,S_0,d_0)\longrightarrow 
\overline{\Cal{M}_X(R_0,i_0,S_0,d_0)}
$$
is continuous.
\endproclaim

Thus, if $(\Mbar_k,g_k), \, (\Mbar,g)\in \Cal{M}_X(R_0,i_0,S_0,d_0)$
and the boundary spectral data of $(\Mbar_k,g_k)$ tend to the boundary
spectral data of $(\Mbar,g)$ in $\Cal{B}$, then, for large $k$,
$\Mbar_k$ are diffeomorphic to $\Mbar$ and $g_k \rightarrow g$ in $C^r,$ 
for all $ r<2$.  


$$\text{}$$
{\bf 5. Auxiliary regularity results}
\newline {}\newline

Here we establish a number of elliptic regularity results, needed
in the analysis in \S{2}, which we did not find in the literature.
These results tend to be variants of known results, but they differ
in various key respects.  Sometimes it is in the category of function 
space involved, e.g., coefficients of the PDE in a non-standard space,
which nevertheless arose naturally in the Ricci equation analysis.
In some cases we can get away with a short argument based on 
standard results, while in other cases we need to do more work.

In \S{5.1} we establish a local regularity result for an elliptic PDE
with coeficients simultaneously satisfying a H{\"o}lder condition and a 
Besov condition.  In \S{5.2} we establish estimates for the Dirichlet 
problem when the coefficients simultaneously satisfy a H{\"o}lder
condition and a Sobolev condition.  In \S\S{5.3--5.5} we obtain estimates
on weak solutions to Neumann boundary problems, with rough coefficients.

$\text{}$

\heading
\S{5.1}: Local Besov regularity
\endheading

$\text{}$

The following result establishes (2.2.8).

\proclaim{Proposition 5.1.1} Assume $u\in H^{1,2}(\Cal{O})$ solves the 
elliptic PDE
$$
\pa_j a^{jk}\pa_k u=0\ \text{ on }\ \Cal{O}.
\tag{5.1.1}
$$
Assume
$$
a^{jk}\in C^r\cap B^s_{p,p}
\tag{5.1.2}
$$
with $r,s\in (0,1),\ p\in (1,\infty)$.  (One should assume $r\le s$.)
Then, locally,
$$
u\in C^{1+r}\cap B^{1+s}_{p,p}.
\tag{5.1.3}
$$
\endproclaim
\demo{Proof} That $u\in C^{1+r}$ is well known; we show $u\in B^{1+s}_{p,p}$.
The proof is like that of Proposition 9.4 in Chapter III of [T2].  
In particular, we use paraproducts, operators of the form $T_au$,
a tool in nonlinear PDE introduced by J.-M.~Bony.  A sketch of the behavior
of paraproducts can be found in Chapter II of [T2].

Set $L^\#=\pa_j T_{a^{jk}}\pa_k$ and write (5.1.1) as
$$
L^\#u=-\pa_j[T_{\pa_ku}a^{jk}+R(a^{jk},\pa_ku)].
\tag{5.1.4}
$$
Here $L^\#\in OP\Cal{B}S^2_{1,1}$ is elliptic and, since $a^{jk}\in C^r$,
by Proposition 6.1 in Chapter I of [T2] we have $E\in OPS^{-2}_{1,1}$ such 
that $EL^\#=I+F,\ F\in OPS^{-r}_{1,1}$.  Then we have
$$
u=-E\pa_j[T_{\pa_ku}a^{jk}+R(a^{jk},\pa_ku)]-Fu.
\tag{5.1.5}
$$
Now
$$
\aligned
u\in C^{1+r}&\Rightarrow T_{\pa_ku},R_{\pa_ku}\in OPS^0_{1,1} \\
&\Rightarrow T_{\pa_ku}a^{jk}+R(a^{jk},\pa_ku)\in B^s_{p,p} \\
&\Rightarrow E\pa_j[T_{\pa_ku}a^{jk}+R(a^{jk},\pa_ku)]\in B^{1+s}_{p,p},
\endaligned
\tag{5.1.6}
$$
the second implication using the hypothesis that $a^{jk}\in B^s_{p,p}$.
Thus (5.1.5) gives
$$
u=-Fu\ \text{ mod }\ B^{1+s}_{p,p},
\tag{5.1.7}
$$
and since $F\in OPS^{-r}_{1,1}$ an iteration from $u\in C^{1+r}$ readily 
yields $u\in B^{1+s}_{p,p}$.
\enddemo

$\text{}$

\heading
\S{5.2}: $L^p$-Sobolev estimates for the Dirichlet problem
\endheading

$\text{}$

Let $\Ombar$ be a smooth, $n$-dimensional, manifold with boundary, 
with metric tensor
$$
g_{jk}\in C^s(\Ombar)\cap H^{1,p}(\Omega),\quad s\in (0,1),\ p\in (1,\infty).
\tag{5.2.1}
$$
We aim to prove the following regularity result.

\proclaim{Proposition 5.2.1} Let $\Cal{O}\subset\bW$ be open.  
Assume $u\in \Lip(\Ombar)$ satisfies
$$
\Delta u=f\in L^p(\Omega),\quad u\bigr|_{\Cal{O}}=0.
\tag{5.2.2}
$$
Also assume $p\ge n/(n-1)$.
Given $z\in \Cal{O}$, there exists a neighborhood $\Ubar$ of $z$ in $\Ombar$
such that
$$
u\in H^{2,p}(U).
\tag{5.2.3}
$$
\endproclaim
\demo{Proof} Take $\varphi\in C_0^\infty(\Ombar)$ such that $\varphi\equiv 1$
near $z$ and $\varphi\equiv 0$ on a neighborhood in $\Ombar$ of $\bW\setminus
\Cal{O}$.  Suppose $\varphi$ is supported in a coordinate patch.  With $a^{jk}
=g^{1/2}g^{jk}$, we have
$$
\pa_ja^{jk}\pa_k(\varphi u)=\varphi g^{1/2}f
+u (\pa_j a^{jk}\pa_k\varphi)
+2a^{jk}(\pa_k\varphi)(\pa_ju),
\tag{5.2.4}
$$
and if (5.2.1) holds and $u\in \Lip(\Ombar)$, each term 
on the right side of (5.2.4) belongs to $L^p(\Omega)$.  
This reduces us to the case $\Omega=\RR^n_+$,
with $u|_{\bW}=0$ and $u$ having compact support in $\Ombar$, where we relabel
$\varphi u$ as $u$.  (This part of the argument still works even if we weaken 
the hypothesis $u\in\Lip(\Ombar)$
to $u\in L^\infty(\Omega)\cap H^{1,p}(\Omega)$.)

Computing formally, we have for $1\le\ell\le n-1$ that $u_\ell=\pa_\ell u$
satisfies
$$
\pa_j a^{jk}\pa_k u_\ell=\pa_\ell(g^{1/2}f)-\pa_j\bigl((\pa_\ell a^{jk})
\pa_ku\bigr).
\tag{5.2.5}
$$
We have $g^{1/2}f\in L^p(\Omega)$ and $(\pa_\ell a^{jk})(\pa_ku)\in 
L^p(\Omega)$, under our hypotheses, so the right side of (5.2.5) belongs to
$H^{-1,p}(\Omega)$, and $u_\ell|_{\bW}=0$. 

We claim this implies
$$
\pa_\ell u\in H^{1,p}(\Omega),\quad 1\le \ell\le n-1.
\tag{5.2.6}
$$
This follows from Theorem 5.5.5' of [Mo2], as long as $p\ge n/(n-1)$.
Granted this,
we have $\pa_j\pa_ku\in L^p(\Omega)$ for all $j,k$ except $j=k=n$, and the
standard trick of using the PDE (5.2.2) to solve for $\pa_n^2u$ yields 
$\pa_n^2u\in L^p(\Omega)$, completing the proof.
\enddemo

Noting some alternative conditions that imply the right sides of (5.2.4) and 
(5.2.5) belong to $L^p(\Omega)$ and $H^{-1,p}(\Omega)$, respectively,
we have the following extension of Proposition 5.2.1.

\proclaim{Proposition 5.2.2}  The conclusion (5.2.3) holds for a solution to 
(5.2.2) provided
$$
g_{jk}\in C^s(\Ombar)\cap H^{1,a}(\Omega),\quad u\in L^\infty(\Omega)\cap
H^{1,b}(\Omega),
\tag{5.2.7}
$$
with $s\in (0,1),\ a,b\in (1,\infty)$, and
$$
\frac{1}{a}+\frac{1}{b}\le \frac{1}{p}.
\tag{5.2.8}
$$
\endproclaim

Returning to the setting of Proposition 5.2.1, we note the following simple
corollary, which is directly applicable to establish (2.2.9).

\proclaim{Corollary 5.2.3} The conclusion of Proposition 5.2.1 remains valid 
if (5.2.2) is generalized to
$$
\Delta u=f\in L^p(\Omega),\quad u\bigr|_{\Cal{O}}=g\in B^{2-1/p}_{p,p}
(\Cal{O})\cap \Lip(\Cal{O}).
\tag{5.2.9}
$$
\endproclaim
\demo{Proof} After perhaps shrinking $\Cal{O}$, we can assume $g=G|_{\Cal{O}}$
with $G\in H^{2,p}(\Omega)\cap \Lip(\Ombar)$.  Then $u-G$ solves
$$
\Delta(u-G)=\tilde{f}\in L^p(\Omega),\quad u-G\bigr|_{\Cal{O}}=0,
\tag{5.2.10}
$$
and Proposition 5.2.1 applies to $u-G$.
\enddemo

$\text{}$

\heading
\S{5.3}: Regularity for weak solutions to the Neumann problem
\endheading

$\text{}$

Let $M$ be a smooth, compact, connected manifold, of dimension $n$.
Assume $M$ has a Riemannian metric tensor that is H{\"o}lder continuous,
of class $C^r$, for some $r\in (0,1)$.  Let $\Omega\subset M$ be a connected
open set, with boundary of class $C^{1+r}$.  Actually we can assume
$\bW$ is smooth, since a $C^{1+r}$ diffeomorphism can smooth out $\bW$
while producing a new metric tensor of class $C^r$.  Let $\Delta$ denote
the Laplace operator on $M$.  Assume $V\in L^\infty(M),\ V\ge 0$ on $M$, and 
$V>0$ on a set of positive measure.  Consider $L=\Delta-V$.  

A weak solution to the Neumann problem
$$
Lu=f\ \text{ on }\ \Omega,\quad Nu=g\ \text{ on }\ \bW
\tag{5.3.1}
$$
is an element $u\in H^1(\Omega)$ satisfying
$$
\int\limits_\Omega \langle du,d\psi\rangle\, dV=
-\int\limits_\Omega (Vu+f)\psi\, dV-\int\limits_{\bW}g\psi\, dS,
\tag{5.3.2}
$$
for all $\psi\in H^1(\Omega)$.  Here the volume element $dV$ on $\Omega$
and the area element $dS$ on $\bW$ are determined by the Riemannian metric
tensor on $M$, as is the inner product $\langle \xi,\eta\rangle$ of 
1-forms.  We aim to prove the following:

\proclaim{Theorem 5.3.1} Given $s\in (0,r),\ u\in H^1(\Omega)$ satisfying 
(5.3.2), $p\ge n/(1-s)$,
$$
f\in L^p(\Omega),\ g\in C^s(\bW)\Longrightarrow u\in C^{1+s}(\Ombar).
\tag{5.3.3}
$$
\endproclaim

Our first reduction is to show that it suffices to take $f=0$.  Indeed, 
extending $f$ by $0$ on $M\setminus\Omega$ we can solve $Lv=f$ on $M$.  By 
Proposition 2.3 of [MT], we have $v\in C^{1+s}(\Ombar)$, provided 
$L^p(M)\subset C^{-1+s}_*(M)$, which holds as long as $H^{1-s,p}(M)\subset
C^0_*(M)$, i.e., as long as $p(1-s)\ge n$.  Then $v|_\Omega$ satisfies
$$
Lv=f\ \text{ on }\ \Omega,\quad Nv=g_0\in C^s(\bW),
\tag{5.3.4}
$$
and it suffices to show that $w=u-v$ belongs to $C^{1+s}(\Ombar)$.

Our next step is to look at
$$
Lw=0\ \text{ on }\ \Omega,\quad Nw=g_1\in C^s(\bW),
\tag{5.3.5}
$$
where $g_1=g-g_0$, and produce a solution $w\in C^{1+s}(\Ombar)$.  
(If $V\equiv 0$ on $\Omega$, assume $\int_{\bW}g_1\, dS=0$.)
Producing such a solution to (5.3.5) will prove Theorm 5.3.1, since a
solution $w\in H^1(\Omega)$ is unique (up to an additive constant if
$V\equiv 0$ on $\Omega$).  To produce such a solution, we use the method of 
layer potentials.

Thus let $E(x,y)$ be the integral kernel of $L^{-1}:H^{-1}(M)\rightarrow
H^1(M)$, and define the single layer potential
$$
\Cal{S}h(x)=\int\limits_{\bW} E(x,y)h(y)\, dS(y).
\tag{5.3.6}
$$
To proceed we need an analysis of $E(x,y)$.
It is elementary to show that $h\in L^1(\bW)\Rightarrow \Cal{S}h\in 
C^{1+r}_{\text{loc}}(M\setminus\bW)$.  In fact, given compact $\Sigma\subset
\bW$,
$$
h\in L^1(\bW),\ \text{supp}\, h\subset\Sigma\Longrightarrow
\Cal{S}h\in C^{1+r}_{\text{loc}}(M\setminus\Sigma).
\tag{5.3.7}
$$
This permits us to use partition of unity arguments and localize the 
study of $\Cal{S}h$ to coordinate patches.  We can choose local coordinates
such that $\bW$ is given by $\{x:x^n=0\}$.

As in [MT], we can write
$$
E(x,y)\sqrt{g(y)}=e_0(x-y,y)+e_1(x,y),
\tag{5.3.8}
$$
where (if $n\ge 3$)
$$
e_0(x-y,y)=C_n\Bigl(\sum g_{jk}(y)(x_j
-y_j)(x_k-y_k)\Bigr)^{-(n-2)/2}.
\tag{5.3.9}
$$
As shown in Theorem 2.6 of [MT], we have, for each $\ep>0$,
$$
\aligned
|e_1(x,y)|&\le C_\ep |x-y|^{-(n-2-r+\ep)}, \\
|\nabla_x e_1(x,y)|&\le C_\ep |x-y|^{-(n-1-r+\ep)}.
\endaligned
\tag{5.3.10}
$$
Also (2.67) of [MT] implies (for $0<s<r$)
$$
|\nabla_x e_1(x_1,y)-\nabla_x e_1(x_2,y)|\le C_\ep
|x_1-x_2|^s\, |x_1-y|^{-(n-1-r+s+\ep)},
\tag{5.3.11}
$$
provided $|x_1-x_2|\le (1/2)|x_1-y|$.

Given $h$ supported in $\Sigma$, the intersection of $\bW$ with a 
coordinate patch, we analyze $\Cal{S}h$ as a sum of two pieces, $\Cal{S}h=
\Cal{S}_0h+\Cal{S}_1h$, where
$$
\aligned
\Cal{S}_0h(x)&=\int\limits_{\bW} e_0(x-y,y)g(y)^{-1/2} h(y)\, dS(y), \\
\Cal{S}_1h(x)&=\int\limits_{\bW} e_1(x,y)g(y)^{-1/2} h(y)\, dS(y).
\endaligned
\tag{5.3.12}
$$

\proclaim{Lemma 5.3.2} Given $s\in (0,r)$, we have
$$
\Cal{S}_1:L^\infty(\Sigma)\longrightarrow C^{1+s}(M).
\tag{5.3.13}
$$
\endproclaim
\demo{Proof} We need to show that, for $x_j\in M,\ h\in L^\infty(\bW)$,
supported in $\Sigma$,
$$
|\nabla\Cal{S}_1h(x_1)-\nabla\Cal{S}_1h(x_2)|\le C|x_1-x_2|^s\, 
\|h\|_{L^\infty}.
\tag{5.3.14}
$$
There are two cases to consider.
\newline {}\newline
{\smc Case I}.  $|x_1-x_2|\le (1/2)\text{dist}(x_1,\Sigma)$.
\newline
Then use (5.3.11) to get (5.3.14).
\newline {}\newline
{\smc Case II}.  $|x_1-x_2|\ge (1/2)\text{dist}(x_1,\Sigma)$.
\newline
Set $\Cal{O}=\{y\in\Sigma:|x_1-y|\le 4|x_1-x_2|\}$.  Use (5.3.10) for $y\in
\Cal{O}$ to analyze separately $\nabla_x \int_{\Cal{O}} e_0(x_j-y,y)
g(y)^{-1/2}h(y)\, dS(y)$, and use (5.3.11) for $y\in\Sigma\setminus\Cal{O}$
to complete the analysis of the left side of (5.3.14).
\enddemo

\proclaim{Lemma 5.3.3} Given $s\in (0,r)$, we have
$$
\Cal{S}_0:C^s(\bW)\longrightarrow C^{1+s}(\Ombar).
\tag{5.3.15}
$$
\endproclaim
\demo{Proof} This is a standard layer potential estimate.  One has
$$
\|\nabla \Cal{S}_0h\|_{L^\infty(\Omega)}\le C\|h\|_{C^s(\bW)},
\tag{5.3.16}
$$
and
$$
|\pa^2 \Cal{S}_0h(x',x_n)|\le C|x_n|^{s-1}\, \|h\|_{C^s(\bW)},
\tag{5.3.17}
$$
which gives 
$$
|\nabla \Cal{S}_0h(x_1)-\nabla\Cal{S}_0h(x_2)|\le C|x_1-x_2|^s\, 
\|h\|_{C^s(\bW)},\quad x_j\in\Ombar;
\tag{5.3.18}
$$
compare Proposition 8.7 in Chapter 13 of [T1].
\enddemo

The next step in solving (5.3.5) in the form $w=\Cal{S}h$ is to analyze
$N\Cal{S}h|_{\bW}$.  We have the standard formula
$$
N\Cal{S}h\bigr|_{\bW}=\Bigl(-\frac{1}{2}I+K^*\Bigr)h,
\tag{5.3.19}
$$
with
$$
K^*h(x)=\text{PV} \int\limits_{\bW} N_x E(x,y)h(y)\, dS(y),\quad x\in \bW.
\tag{5.3.20}
$$
Compare (2.81) of [MT] for such a formula in the more general context
of a Lipschitz boundary $\bW$.
We can break $K^*$ into two pieces, $K^*=K_0^*+K_1^*$, using (5.3.8):
$$
K^*_0h(x)=\text{PV}\int\limits_{\bW} N_x e_0(x-y,y)g(y)^{-1/2}h(y)\, dS(y),
\tag{5.3.21}
$$
and
$$
K^*_1h(x)=\int\limits_{\bW} N_x e_1(x,y)g(y)^{-1/2} h(y)\, dS(y).
\tag{5.3.22}
$$
By Lemma 5.3.2 we have $K_1^*:L^\infty(\bW)\rightarrow C^s(\bW)$.  As for
$K^*_0$, it is a pseudodifferential operator with double symbol, of the 
sort studied in Chapter I, \S{9} of [T2], with symbol
$$
a(x,y,\xi)\in C^rS^0_{1,0}.
\tag{5.3.23}
$$
Furthermore, the appearance of $N_x$ in (5.3.21) produces the following
important cancellation effect:
$$
a(x,x,\xi)=0.
\tag{5.3.24}
$$
(This is part of what makes analysis on domains with $C^{1+r}$ boundary 
easier than analysis on Lipschitz domains.)  Hence, by Proposition 9.15
of [T2], Chapter I, we have $K^*_0:L^\infty(\bW)\rightarrow C^s(\bW)$, 
for all $s<r$.  In summary,
$$
K^*:L^\infty(\bW)\longrightarrow C^s(\bW),\quad \forall\ s<r,
\tag{5.3.25}
$$
and hence
$$
K^*:C^s(\bW)\longrightarrow C^s(\bW)\ \text{ is compact, for all }\
s\in (0,r).
\tag{5.3.26}
$$
Thus
$$
-\frac{1}{2}I+K^*:C^s(\bW)\longrightarrow C^s(\bW)\ \text{ is Fredholm, of
index }\ 0.
\tag{5.3.27}
$$
Using this we can prove the following.

\proclaim{Proposition 5.3.4} Let $g_1\in C^s(\bW),\ 0<s<r$.  
If $V>0$ somewhere 
(i.e., on a set of positive measure) on $\Omega$, 
then (5.3.5) has a unique solution $w=\Cal{S}h\in C^{1+s}(\Ombar)$, where
$$
\Bigl(-\frac{1}{2}I+K^*\Bigr)h=g_1,\quad h\in C^s(\bW).
\tag{5.3.28}
$$
If $V\equiv 0$ on $\Omega$, then (5.3.5) has a solution $w\in C^{1+s}(\Ombar)$
if and only if $\int_{\bW}g_1\, dS=0$, and such $w$ is unique up to an
additive constant.
\endproclaim
\demo{Proof} This is a standard argument in layer potential theory.
One shows that, if $V>0$ somewhere on $\Omega$, then $-(1/2)I+K^*$ is
injective on $C^s(\bW)$, hence, by (5.3.27), bijective.  If $V\equiv 0$ on 
$\Omega$, then $-(1/2)I+K^*$ has a one-dimensional kernel in $C^s(\bW)$,
and the constant function $1$ annihilates its range.  Compare the treatment
of Theorem 3.4 in [MT], carried out in the more general context of a
Lipschitz boundary.  Once one solves (5.3.28) for $h\in C^s(\bW)$, the fact 
that $w=\Cal{S}h$ solves (5.3.5) and belongs to $C^{1+s}(\Ombar)$ follows 
from the previous analysis.
\enddemo

With this result, the proof of Theorem 5.3.1 is complete.

$\text{}$

\heading
\S{5.4}: Local regularity for the Neumann problem
\endheading

$\text{}$

It is useful to strengthen the global regularity result of \S{5.3} 
to a local regularity result for a weak solution to the
Neumann problem.  With $M, \Omega$, and $L$ as in \S{5.3}, let $\Cal{O}$
be an open subset of $\bW$ and suppose $u$ is a weak solution to
$$
Lu=f\ \text{ on }\ \Omega,\quad Nu\bigr|_{\Cal{O}}=g.
\tag{5.4.1}
$$
That is to say, we assume $u\in H^1(\Omega)$ and that (5.3.2) holds for
all $\psi\in H^1(\Omega)$ such that $\psi|_{\bW}$ vanishes on a neighborhood
of $\bW\setminus \Cal{O}$ in $\bW$.  We prove the following.

\proclaim{Theorem 5.4.1} Let $u$ be a weak solution to (5.4.1).  Assume $u\in 
H^{1,q}(\Omega)$ with $q\ge 2$ and either $q>1/r$ or $q\ge p$. 
(Recall $g_{jk}\in C^r$.)
As in Theorem 5.3.1, assume $0<s<r<1$ and $p\ge n/(1-s)$.  Then
$$
f\in L^p(\Omega),\ g\in C^s(\Cal{O})\Longrightarrow 
u\in C^{1+s}(\Omega\cup\Cal{O}).
\tag{5.4.2}
$$
\endproclaim
\demo{Proof} Take $p_0\in\Cal{O}$ and pick $\varphi\in C_0^\infty(M)$, 
equal to $1$ on a neighborhood of $p_0$, such that $\varphi=0$ on a 
neighborhood of $\bW\setminus\Cal{O}$ in $M$.  Let $v=\varphi u|_{\Omega}$,
so $v\in H^{1,q}(\Omega)$.  We seek to establish extra regularity of $v$.
It is readily verified that $v$ is a (global) weak solution of
$$
Lv=\tilde{f}\ \text{ on }\ \Omega,\quad Nv=\tilde{g}\ \text{ on }\ \bW,
\tag{5.4.3}
$$
with 
$$
\aligned
\tilde{f}&=\varphi f+2\langle d\varphi,du\rangle+u\Delta \varphi, \\
\tilde{g}&=\varphi g+u(N\varphi)\bigr|_{\bW}.
\endaligned
\tag{5.4.4}
$$
We have
$$
\tilde{f}\in L^{p\wedge q}(\Omega),\quad \tilde{g}\in C^s(\bW)+
C^r(\bW)\cdot B^{1-1/q}_{q,q}(\bW).
\tag{5.4.5}
$$
Extend $\tilde{f}$ by $0$ on $M\setminus \Omega$ and set $v_1=L^{-1}
\tilde{f}$.  By (2.16) of [MT], we have
$$
L^{p\wedge q}(M)\subset H^{\rho-1,\tilde{p}}(M)\Longrightarrow
v_1\in H^{\rho+1,\tilde{p}}(M),\quad \rho=r-\ep.
\tag{5.4.6}
$$
In fact, if $q\ge p$, we can use Proposition 2.3 of [MT] as in the beginning
of the proof of Theorem 5.3.1, and say $v_1\in C^{1+s}(M)$.  
If $q<p$, then we can take $\tilde{p}>q$ in (5.4.6).

Now $v_2=v-v_1|_{\Omega}$ is a weak solution to
$$
Lv_2=0\ \text{ on }\ \Omega,\quad Nv_2=g_2=\tilde{g}-Nv_1,
\tag{5.4.7}
$$
and
$$
Nv_1\in C^r(\bW)\cdot B^{\rho-1/\tilde{p}}_{\tilde{p},\tilde{p}}(\bW),
\tag{5.4.8}
$$
if $q<p$, and in $C^s(\bW)$ if $q\ge p$.  Note that $\tilde{p}\rho>q(r-\ep)$
can be assumed to be $>1$ by the hypothesis $rq>1$, if we take $\ep>0$
small enough.  At this point we have $g_2\in L^{\tilde{q}_2}(\bW)$, with
$\tilde{q}_2>q$.  From this we can deduce
$$
v_2\in H^{1,\tilde{q}_2}(\Omega),\quad \tilde{q}_2>q.
\tag{5.4.9}
$$

The result (5.4.9) is a special case of much stronger known results; let us
sketch the proof.  First, parallel to (5.3.26), we have
$$
K^*:L^p(\bW)\longrightarrow L^p(\bW)\ \text{ is compact, for all }
p\in (1,\infty).
\tag{5.4.10}
$$
Indeed, the compactness of $K_0^*$ follows from Proposition 9.5 in Chapter I
of [T2], together with (5.3.23)--(5.3.24), and the compactness of $K^*_1$
follows from the estimates in (5.3.10).  Having (5.4.10), we see that
$$
-\frac{1}{2}I+K^*:L^p(\bW)\longrightarrow L^p(\bW)\ \text{ is Fredholm,
of index }0,
\tag{5.4.11}
$$
for each $p\in (1,\infty)$.  Then an argument as in the proof of Proposition 
5.3.4 yields
$$
v_2=\Cal{S}h_2,\quad h_2=\Bigl(-\frac{1}{2}I+K^*\Bigr)^{-1}g_2
\in L^{\tilde{q}_2}(\bW),
\tag{5.4.12}
$$
up to an additive constant if $V\equiv 0$ on $\Omega$.  Now we have the
non-tangential maximal function estimate
$$
\|(\nabla \Cal{S}h_2)^*\|_{L^{\tilde{q}_2}(\bW)}\le C_{\tilde{q}_2}
\|h_2\|_{L^{\tilde{q}_2}(\bW)},\quad 1<\tilde{q}_2<\infty,
\tag{5.4.13}
$$
which is stronger than (5.4.9).  The estimate (5.4.13), in the setting of
(5.4.6) for a $C^r$-metric tensor, is proven in (2.77) of [MT], in the
more general context of a Lipschitz domain $\Omega$.  

Now (5.4.6) and (5.4.9) together give
$$
v=v_1+v_2\in H^{1,q_2}(\Omega),\quad q_2>q.
\tag{5.4.14}
$$
Now we go back to $u$.  After shrinking $\Omega$ to a smaller neighborhood
of $p_0$, we can replace our hypothesis $u\in H^{1,q}(\Omega)$ by 
$u\in H^{1,q_2}(\Omega)$.  Iterating this argument yields $u\in H^{1,q_\nu}
(\Omega)$ with $q<q_2<q_3<\cdots$.  After a finite number of iterations we
reach a point where Theorem 5.3.1 is applicable to $v=\varphi u$, and
Theorem 5.4.1 is proven.
\enddemo

$\text{}$

\heading
\S{5.5}: Neumann data in $L^s$
\endheading

$\text{}$

We produce more regularity results on weak solutions to
$$
Lu=f\ \text{ on }\ \Omega,\quad Nu=g\ \text{ on }\ \bW,
\tag{5.5.1}
$$
starting with the following.  As before, assume $g_{jk}\in C^r$ for some
$r>0$.

\proclaim{Proposition 5.5.1} Let $u\in H^{1,2}(\Omega)$ be a weak solution 
to (5.5.1).  Assume
$$
f\in L^p(\Omega),\quad g\in L^s(\bW),\quad p>n,\ 1<s<\infty.
\tag{5.5.2}
$$
Then we have
$$
(\nabla u)^*\in L^s(\bW),\quad u\in H^{1,ns/(n-1)}(\Omega).
\tag{5.5.3}
$$
\endproclaim
\demo{Proof} Since $p>n$, we can solve $Lv=f$ (with $f$ extended by $0$)
on a neighborhood $\Cal{O}$ of $\Ombar$, with $v\in C^{1+\sigma}(\Cal{O}),\
\sigma>0$.  This reduces our consideration to the case $f=0$ in (5.5.2).
Then $u$ is given by
$$
u=\Cal{S}h,
\tag{5.5.4}
$$
with
$$
\Bigl(-\frac{1}{2}I+K^*\Bigr)h=g,\quad h\in L^s(\bW).
\tag{5.5.5}
$$
The fact that
$$
h\in L^s(\bW)\Longrightarrow (\nabla \Cal{S}h)^*\in L^s(\bW)
\tag{5.5.6}
$$
is established in the context of a H{\"o}lder continuous metric tensor
(and in the more general context of a Lipschitz domain) in [MT].  

The last part of (5.5.3) follows from the mapping property
$$
\Cal{S}:L^s(\bW)\longrightarrow H^{1,ns/(n-1)}(\Omega).
\tag{5.5.7}
$$
This is demonstrated, in the more general context of a Lipschitz domain,
in [MT2].  We describe here the basic structure of the argument.
We write $\Cal{S}=\Cal{S}_0+\Cal{S}_1$, as in (5.3.12).  
Harmonic analysis techniques are brought to bear to establish
$$
\Cal{S}_0:L^s(\bW)\longrightarrow H^{1,ns/(n-1)}(\Omega).
\tag{5.5.8}
$$
(Such a result is easier for a smooth domain than for a Lipschitz domain.)
As for $\Cal{S}_1$, we already have in (5.3.13) that
$$
\Cal{S}_1:L^\infty(\bW)\longrightarrow C^{1+\sigma}(\Ombar),
\tag{5.5.9}
$$
for some $\sigma>0$.  Meanwhile the estimate (5.3.10) on $\nabla_xe_1(x,y)$ is
more than adequate to give
$$
\Cal{S}_1:L^1(\bW)\longrightarrow H^{1,n/(n-1)}(\Omega),
\tag{5.5.10}
$$
and then interpolation gives more than
$$
\Cal{S}_1:L^s(\bW)\longrightarrow H^{1,ns/(n-1)}(\Omega).
\tag{5.5.11}
$$

We now establish a useful local regularity result.

\proclaim{Proposition 5.5.2} Let $\Cal{O}$ be an open subset of $\bW$ and 
assume $u$ is a weak solution to
$$
Lu=f\ \text{ on }\ \Omega,\quad Nu\bigr|_{\Cal{O}}=g.
\tag{5.5.12}
$$
Take $p>n,\ s\in (1,\infty)$, and assume 
$$
u\in H^{1,p}(\Omega),\quad f\in L^p(\Omega),\quad g\in L^s(\Cal{O}).
\tag{5.5.13}
$$
Then each $p\in \Cal{O}$ has a neighborhood $\Ubar$ in $\Ombar$ such that
$$
u\in H^{1,ns/(n-1)}(U).
\tag{5.5.14}
$$
\endproclaim
\demo{Proof} As in the proof of Theorem 5.4.1, we consider $v=\varphi u$, 
which is a global weak solution of (5.4.3), with $\tilde{f}, \tilde{g}$ 
given by (5.4.4).  Under our current hypotheses we have
$$
\tilde{f}\in L^p(\Omega),\quad \tilde{g}\in L^s(\bW),
\tag{5.5.15}
$$
so Proposition 5.5.1 gives $v\in H^{1,ns/(n-1)}(\Omega)$, and Proposition 
5.5.2 is proven.
\enddemo

$$\text{}$$
{\bf References}

\roster
\item"[Al]" G.~Alessandrini, Stable determination of conductivity by boundary 
measurements, Applic. Anal. 27 (1988), 153--172.
\item"[An1]" M.~Anderson, Convergence and rigidity of manifolds under
Ricci curvature bounds, Invent. Math. 102 (1990), 429--445.
\item"[An2]" M.~Anderson, Boundary regularity, uniqueness, and non-uniqueness
for AH Einstein metrics on 4-manifolds, Preprint, 2001.
\item"[Be1]" M.~Belishev, An approach to multidimensional inverse problems
(Russian), Dokl. Akad. Nauk. SSSR, 297 (1987), 524--527.
\item"[BK1]" M.~Belishev and Y.~Kurylev, To the reconstruction of a 
Riemannian manifold via its spectral data (BC-method), Comm. PDE 17 (1992),
767--804.
\item"[Bz]" Y.~Berezanskii, The uniqueness theorem in the inverse problem
of spectral analysis for the Schr\"odinger equation (Russian), 
Trudy Moskov. Matem. Obsch., 7 (1958), 1--62.
\item"[Ch]" J.~Cheeger, Finiteness theorems for Riemannian manifolds,
Amer. J. Math. 92 (1970), 61--75.
\item"[CG]" J.~Cheeger and D.~Gromoll, The splitting theorem for manifolds
with non-negative Ricci curvature, J. Diff. Geom. 6 (1971), 119--128.
\item"[CW]" Y.-Z.~Chen and L.-C.~Wu, Second Order Elliptic Equations and
Elliptic Systems, Transl. Math. Monogr., No.~174, AMS, Providence, R.I.,
1998.
\item"[DTK]" D.~De Turck and J.~Kazdan, Some regularity theorems in 
Riemannian geometry, Ann. Scient. Ecole Norm. Sup. Paris, 14 (1981),
249--260.
\item"[Ge]" I.~Gel'fand, Some aspects of functional analysis and algebra,
Proc. ICM 1 (1954), 253--277.
\item"[GT]" D.~Gilbarg and N.~Trudinger, Elliptic Partial Differential
Equations of Second Order, Springer-Verlag, New York, 1983.
\item"[GW]" R.~Greene and H.~Wu, Lipschitz convergence of Riemannian 
manifolds, Pacific J. Math. 131 (1988), 119--141.
\item"[Gr]" M.~Gromov, Metric Structures for Riemannian and Non-Riemannian
Spaces (with appendices by M.~Katz, P.~Pansu, and S.~Semmes), Birkhauser,
Boston, 1999.
\item"[Ha]" P.~Hartman, On the local uniqueness of geodesics, Amer. J.
Math. 72 (1950), 723--730.
\item"[HH]" E.~Hebey and M.~Herzlich, Harmonic coordinates, harmonic
radius, and convergence of Riemannian manifolds, Rend. di Mathem. Ser.
VII, 17 (1997), 569--605.
\item"[Ho]" L.~H{\"o}rmander, On the uniqueness of the Cauchy problem 
under partial analyticity assumptions.  In: Geom. Optics and Related
Topics (F.~Colombini and N.~Lerner, eds.), Birkhauser, Boston, 1997.
\item"[KKL]" A.~Katchalov, Y.~Kurylev, and M.~Lassas, Inverse Boundary
Spectral Problems, Chapman-Hall/CRC Press, Boca Raton, 2001.
\item"[KLM]" A.~Katchalov, Y.~Kurylev, M.~Lassas, and N.~Mandache,
Equivalence of time-domain inverse problems and boundary spectral problems,
Preprint, 2002.
\item"[K]" T.~Kato, Perturbation Theory for Linear Operators,
Springer-Varlag, New York, 1966.
\item"[Ka]" A.~Katsuda, BC-method and stability of Gel'fand inverse spectral 
problem, Proc. Conf. ``Spectral and Scattering Theory,'' RIMS, Kyoto,
2001, 24--35.
\item"[K2L]" A.~Katsuda, Y.~Kurylev, and M.~Lassas, Stability in Gelfand
inverse boundary spectral problem, Preprint, 2001.
\item"[Kod]" S.~Kodani, Convergence theory for Riemannian manifolds with
boundary, Compos. Math. 75 (1990), 171--192.
\item"[Ku]" Y.~Kurylev, Multidimensional Gel'fand inverse boundary problem
and boundary distance map.  In: Inverse Problems Related to Geometry
(H.~Soga, ed.), 1--15, Ibaraki Univ. Press, Mito, 1997.
\item"[KuL]" Y.~Kurylev and M.~Lassas, The multidimensional Gel'fand inverse
problem for non-self-adjoint operators, Inverse Problems, 13 (1997), 
1495--1501.
\item"[Mi]" C.~Miranda, Partial Differential Equations of Elliptic Type,
Springer-Verlag, New York, 1970.
\item"[MT]" M.~Mitrea and M.~Taylor, Potential theory on Lipschitz domains
in Riemannian manifolds: H{\"o}lder continuous metric tensors, Comm. PDE
27 (2000), 1487--1536.
\item"[MT2]" M.~Mitrea and M.~Taylor, Local and global regularity on 
Lipschitz domains, Manuscript, in preparation.
\item"[Mo1]" C.~Morrey, Second order elliptic systems of differential 
equations, pp.~101--159 in Ann. of Math. Studies \#33, Princeton Univ. Press,
Princeton, N.J., 1954.
\item"[Mo2]" C.~Morrey, Multiple Integrals in the Calculus of Variations,
Springer-Verlag, New York, 1966.
\item"[NSU]"  A.~Nachman, J.~Sylvester and  G.~Uhlmann, 
An $n-$dimensional Borg--Levinson theorem,
Comm. Math. Phys. 115 (1988),
595-605.
\item"[Nv]"  R.~Novikov, A multidimensional inverse spectral problem 
for the equation $-\Delta \psi + (v(x) -Eu(x)) \, \psi =0$ (Russian),
Funkt. Anal. i Priloz. 22 (1988), 11--22. 
\item"[P]" S.~Peters, Convergence of Riemannian manifolds, Compos. Math.
62 (1987), 3--16.
\item"[Pe]" P.~Petersen, Riemannian Geometry, Springer-Verlag, New York,
1998.
\item"[Sak]" T.~Sakai, On continuity of injectivity radius function,
Math. J. Okayama Univ. 25 (1983), 91--97.
\item"[StU]" P.~Stefanov and G.~Uhlmann, Stability estimates for the
hyperbolic Dirichlet-to-Neumann map in anisotropic media, J. Funct. Anal.
154 (1998), 330--357.
\item"[St]" E.~Stein, Topics in Harmonic Analysis Related to the 
Littlewood-Paley Theory, Princeton Univ. Press, Princeton, N.J., 1970.
\item"[Ta]" D.~Tataru, Unique continuation for solutions to PDE's:
between H{\"o}rmander's theorem and Holmgren's theorem, Comm. PDE 20
(1995), 855--884.
\item"[Ta2]" D.~Tataru, Unique continuation for solutions to PDEs with
partially analytic coefficients, J. Math. Pures Appl. 78 (1999), 505--521.
\item"[T1]" M.~Taylor, Partial Differential Equations, Vols.~1--3, 
Springer-Verlag, New York, 1996.
\item"[T2]" M.~Taylor, Tools for PDE, AMS, Providence, R.I., 2000.
\item"[Tr1]" H.~Triebel, Interpolation Theory, Function Spaces, Differential
Operators, DVW, Berlin, 1978.
\item"[Tr2]" H.~Triebel, Function Spaces, Birkhauser, Boston, 1983.

\endroster

$$\text{}$$
{\smc Michael Anderson}, Mathematics Department, State Univ. of New York, 
Stony Brook, N.Y., 11794, USA
\newline {}\newline
{\smc Atsushi Katsuda}, Mathematics Department, Okayama Univ., 
Tsushima-naka, Okayama, 700-8530, Japan
\newline {}\newline
{\smc Yaroslav Kurylev}, Department of Math. Sciences, Loughborough Univ., 
Loughborough, LE11 3TU, UK
\newline {}\newline
{\smc Matti Lassas}, Rolf Nevanlinna Institute, Univ. of Helsinki, 
FIN-00014, Finland
\newline {}\newline
{\smc Michael Taylor}, Mathematics Department, Univ. of North Carolina, 
Chapel Hill, N.C., 27599, USA

\end